\input amstex
\documentstyle{amsppt}
\magnification = 1200

\hcorrection{.25in}
\def\({\left(}
\def\){\right)}
\def\[{\left[}
\def\]{\right]}

\def\calH{\Cal H}
\def\frakF{\frak F}
\def\frakG{\frak G}
\def\N{\Bbb N}

\def\C{\Bbb C}

\def\Be{1}
\def\BH{2}
\def\Bo{3}
\def\Boo{4}
\def\Booo{5}
\def\Carle{6}
\def\Car{7}
\def\CG{8}
\def\S{9}
\def\Edw{10}
\def\Gne{11}
\def\HR{12}
\def\HW{13}
\def\Heaa{14}
\def\Hea{15}
\def\Heb{16}
\def\Hi{17}
\def\Kah{18}
\def\La{19}
\def\Sch{20}
\def\Y{21}

\NoBlackBoxes
\TagsOnRight
%\refstyle{C}
\topmatter
\title
A Hilbert space of Dirichlet series and systems of dilated functions
in $L^{2}(0,1)$
\endtitle
\rightheadtext{Dirichlet series and systems of dilated functions
in $L^{2}(0,1)$}
\author
H\aa kan Hedenmalm, Peter Lindqvist, and Kristian Seip
\endauthor
\footnotetext"{}"{The work of the first author was supported in part by
the Swedish Natural Science Research Council (NFR), and that of the second and
third authors by the Norwegian Research Council (NFR). While at MSRI, 
Berkeley, the first and third authors were also supported in part by NSF 
grant DMS-9022140.}
%\affil
%\endaffil
\address
\hskip-\parindent
H\aa kan Hedenmalm, Department of Mathematics, Uppsala University, Box 480,
S--75106 Uppsala, Sweden
\endaddress
\email
haakan\@math.uu.se
\endemail
\address
\hskip-\parindent
Peter Lindqvist, Department of Mathematical Sciences, The Norwegian Institute 
of Technology, N--7034 Trondheim, Norway
\endaddress
\email
lqvist\@imf.unit.no
\endemail
\address
\hskip-\parindent
Kristian Seip, Department of Mathematical Sciences, The Norwegian Institute of
Technology, N--7034 Trondheim, Norway
\endaddress
\email
seip\@imf.unit.no
\endemail
%\dedicatory
%\enddedicatory
%\date
%\enddate
%\thanks
%\endthanks
%\translator
%\endtranslator
%\keywords
%\endkeywords
%\subjclass
%\endsubjclass
\abstract
For a function $\varphi$ in $L^2(0,1)$, extended to the whole real line as an
odd periodic function of period 2, we ask when the collection of dilates 
$\varphi(nx)$, $n=1,2,3,\ldots$, constitutes a Riesz basis or a complete 
sequence in $L^2(0,1)$. The problem translates into a question concerning 
multipliers and cyclic vectors in the Hilbert space $\calH$ of Dirichlet 
series $f(s)=\sum_n a_nn^{-s}$, where the coefficients $a_n$ are square 
summable. It proves useful to model $\calH$ as the $H^2$ space of the 
infinite-dimensional polydisk, or, which is the same, the $H^2$ space of the
character space, where a character is a multiplicative homomorphism from the 
positive integers to the unit circle. For given $f$ in $\calH$ and characters
$\chi$, $f_\chi(s)=\sum_na_n\chi(n)n^{-s}$ is a vertical limit function of 
$f$. We study certain probabilistic properties of these vertical limit 
functions.
\endabstract
\endtopmatter

\def\C{\Bbb C}
\def\Cplus{\C_+}
\def\D{\Bbb D}
\def\Q{\Bbb Q}
\def\T{\Bbb T}

\def\E{\frak E}
\def\P{\frak P}
\def\fQ{\frak Q}

\def\calM{\Cal M}
\def\calA{\Cal A}

\def\vare{\varepsilon}

\def\t{t}
\document
\head
1. Introduction
\endhead

\def\Char{\Xi}

The purpose of this paper is twofold. First, we study systems of functions of
the form $\varphi(x),\varphi(2x),\varphi(3x),...$, and second, we consider
the Hardy space $H^{2}$ of the infinite-dimensional polydisk. Building on
ideas of Arne Beurling and Harald Bohr, we find that the two topics are
intimately connected, the common feature being the use of Dirichlet series.

Let  $\varphi\in L^{2}(0,1)$ be given and consider $\varphi$ as defined on
the whole real axis by extending it to an odd periodic function of period 2.
The Riesz-Fischer theorem of Fourier analysis states that for
$\varphi(x)=\sqrt{2}\sin(\pi x)$ the sequence $\varphi(nx)$, $n=1,2,3,...$,
is an orthonormal basis in the Hilbert space $L^{2}(0,1)$. The question raised
in this paper is which functions can take the place of the sine in this
theorem. It is clear that the statement must be weakened, because the only
orthogonal bases are obtained from $\varphi(x)=C\sin(\pi x)$. If we instead
ask for a classification of those $\varphi$ for which the system $\{\varphi
(nx)\}_n$ is a Riesz basis (a basis orthonormal
with respect to an equivalent norm) or of those $\varphi$ for which
the same system is a complete sequence in $L^{2}(0,1)$, we are led to 
profound problems.

The latter of the two problems -- the completeness problem -- was stated by
Beurling in his seminar on harmonic analysis in Uppsala in 1945. A brief note
from this seminar is found in \cite{\bf\Be}. Beurling's note indicates that a
natural way to approach these problems is to associate to the given function
$$\varphi(x)=\sum_{n=1}^{\infty} a_{n}\sqrt{2}\sin(n\pi x)$$
the Dirichlet series
$$S\varphi(s)=\sum_{n=1}^{\infty} a_{n} n^{-s},\tag 1-1$$
and to try to express the Riesz basis and completeness properties in terms of
analytic properties of $S\varphi(s)$. This approach has proved fruitful. We 
have solved completely the Riesz basis problem (Theorems 3.1 and 5.2 below):
{\sl the system $\{\varphi(nx)\}_n$ is a Riesz basis in $L^{2}(0,1)$ if and
only if \thetag{1-1} defines an analytic function bounded away from zero and
infinity in the half-plane $\Re s>0$.} 

A major part of this paper consists of a study of Dirichlet series
of the form \thetag{1-1} with $\sum_n |a_{n}|^{2}<+\infty$. Such Dirichlet
series form, in a natural way, a Hilbert space of analytic functions on the
half-plane $\Re s>\frac12$. From now on, we denote this space by $\calH$.
In Sections 2--4, we seek to reveal the basic analytic features of this
space. A central problem is to characterize the so-called multipliers of 
$\calH$. Theorem 3.1 below states that {\sl the multipliers are precisely 
the bounded analytic functions in the right half-plane $\Re s>0$ which can be 
represented as Dirichlet series.} This result is crucial for the 
characterization of the Riesz bases $\{\varphi(nx)\}_n$ mentioned above. 
Apart for its relevance for the dilation Riesz basis problem, the study of 
the space $\calH$ can also be motivated by the mere fact that its kernel 
function $K_{\calH}(z,w)$ is more or less the Riemann zeta function: 
$K_{\calH}(z,w)=\zeta(z+\bar w)$.

Following an idea of Bohr, we find it both convenient and illuminating to
use the infinite-dimensional polydisk for the study of $\calH$. This leads
us to power series in infinitely many variables, a concept studied already by
Hilbert, and to an identification of $\calH$ as the Hardy space $H^{2}$ of
the infinite-dimensional polydisk. The infinite-dimensional polycircle can be
identified with the group $\Char$ of complex-valued characters $\chi$ on the
positive integers, which satisfy $|\chi(n)|=1$ and $\chi(mn)=\chi(m)\chi(n)$.
The characters constitute the (compact) dual group of the discrete
multiplicative group of positive rationals $\langle{\Bbb Q}_+,\cdot\rangle$.

A function $f\in\calH$ is said to be {\sl cyclic} provided that the collection
of functions $fg$, where $g$ is a finite Dirichlet series, is dense in 
$\calH$. A basic observation is that the system $\{\varphi_n\}_n$ is 
complete if and only if the Dirichlet series $S\varphi$ is cyclic in $\calH$. 
We restate cyclicity in terms of our model of $\calH$ as the $H^2$ space of 
the infinite-dimensional polydisk, and state some conditions for cyclicity. 
There is some overlap between these conditions and results due to Henry 
Helson \cite{\bf\Heb}, who studied cyclicity in a more general context. 

Since the multipliers of $\calH$ extend analytically and boundedly to $\Re
s>0$, but the functions of $\calH$ need only be analytic on $\Re s>\frac12$,
one is led to suspect that, nevertheless, in a sense to be made precise, a
function in $\calH$, picked at random, almost surely extends holomorphically
to $\Re s>0$. Given a Dirichlet series in $\calH$,
$$f(s)=\sum_{n=1}^{\infty}a_{n} n^{-s},$$
we consider, for characters $\chi$,
$$f_\chi(s)=\sum_{n=1}^{\infty}a_{n}\chi(n)n^{-s},$$
which again is an element of $\calH$. The functions $f_\chi$ are the normal
limits of vertical translates of $f$. Almost every one of them has a 
convergent Dirichlet series on $\Re s>0$, and is hence holomorphic there. 
The "almost every" is with respect to the Haar measure on the character group 
$\Char$, which is normalized to have total mass $1$. This result was obtained 
in a more abstract setting by Helson \cite{\bf\Heb}. A curious consequence 
is that a probabilistic ``Riemann hypothesis'' (Corollary 4.8) holds: for 
almost all characters $\chi$, the Dirichlet series $\zeta_\chi(s)=
\sum_{n=1}^\infty\chi(n)n^{-s}$ converges to a zero-free holomorphic function 
in the region $\Re s>\frac12$. The convergence part has the following 
interpretation: for almost all characters $\chi$,
$$\sum_{n=1}^N\chi(n)=O\big(N^{1/2+\vare}\big)\qquad\text{as}\quad
N\to+\infty$$
holds for fixed $\vare>0$. For each $n$, it is fruitful to regard the function
$\chi(n)$ as a stochastic variable, which for $n>1$ is uniformly distributed
on $\T$. The stochastic variables $\chi(p)$, where $p$ runs through the primes,
are then mutually independent, and the mutual dependence of the $\chi(n)$, 
as $n$ runs through the positive integers, is governed by the multiplicative 
rule $\chi(mn)=\chi(m)\chi(n)$. The above estimate of the partial sum 
$\sum_{n:n\le N}\chi(n)$ is what one would expect if all the $\chi(n)$ were 
mutually independent, so although they do satisfy complicated multiplicative 
dependence relations, this is insignificant additively.

The above-mentioned assertions may be compared with Jean-Pierre Kahane's 
paper \cite{\bf\Kah} on random Dirichlet series. Kahane works with functions
$f_\chi(s)$, where $\chi(n)$ is treated as a sequence of independent random
variables, and hence no multiplicativity property holds.

In \cite{\bf\Hea}, Helson
suggests that the classical theory of Dirichlet series be combined with 
modern techniques from harmonic and functional analysis. We hope that the
present paper can inspire work in that direction.

\head
2. The Hilbert spaces ${\Cal H}$ and $H^{2}(\D^{\infty})$
\endhead

\subhead 2.1 The space $\calH$ and preliminaries on Dirichlet series
\endsubhead We are concerned with Dirichlet series of the form
$$ f(s)=\sum^\infty_{n=1}a_nn^{-s}, \tag 2-1$$
where $s=\sigma+it$ is a complex variable.
We recall briefly some classical facts about such series. There are a number
of critical lines or absciss\ae{} connected with \thetag{2-1}. We have {\sl
the abscissa of absolute convergence} $\sigma_{a}$ and {\sl the abscissa of
ordinary convergence} $\sigma_{c}$. These numbers are such that the series
converges in the prescribed sense to the right but not to the left of the
abscissa in question. We also have {\sl the abscissa of uniform convergence}
$\sigma_{u}$, defined as the infimum of those $\sigma_{0}$ for which the
series converges uniformly in the half-plane $\Re s>\sigma_{0}$. We have
trivially $-\infty\leq\sigma_{c}\leq\sigma_{u}\leq\sigma_{a}\leq+\infty$ and
$\sigma_{a}-\sigma_{c}\leq 1$ if anyone of the absciss\ae{} is finite. A
theorem of Bohr \cite{\bf\Boo} says that $\sigma_{a}-\sigma_{u}\leq\frac12$,
and this inequality is sharp, as was shown by Bohnenblust and Hille
\cite{\bf\BH}.

When we need to indicate which function we are taking the absciss\ae{} of,
we do this by putting the function in parenthesis; for instance, we would
write $\sigma_c(f)$ for the abscissa of convergence. In terms of the partial 
sums of the coefficients, $S_N=\sum_{n:n\le N} a_n$, the abscissa of 
convergence of the Dirichlet series \thetag{2-1} can be expressed as follows: 
unless $\sigma_c(f)$ is negative, $\sigma_c(f)$ equals the infimum of all 
positive real numbers $\alpha$ for which $S_N=O(N^\alpha)$ as $N\to+\infty$.

Of relevance to us is also {\sl the abscissa of regularity and boundedness}
$\sigma_{b}$, which is the infimum of those $\sigma_{0}$ for which the
function defined by \thetag{2-1} (possibly by analytic continuation from a
smaller half-plane) is analytic and bounded for $\Re s>\sigma_{0}$.
We shall need the following result of Bohr \cite{\bf\Bo}.

\proclaim{Lemma 2.1} {\rm(Bohr's theorem)} $\sigma_{u}=\sigma_{b}$.
\endproclaim

For a more complete account of the basic facts on the convergence of Dirichlet
series, the reader is referred to \cite{\bf\HR} and Bohr's thesis 
\cite{\bf\Booo}.

We will assume that the function $f$ given by the Dirichlet series
\thetag{2-1} belongs to $\Cal H$, that is,
$$\sum_{n=1}^{\infty}|a_{n}|^{2}<+\infty.\tag 2-2$$
The Cauchy-Schwarz inequality yields
$$\Big(\sum^\infty_{n=1}\big|a_nn^{-s}\big|\Big)^2\leq\sum^\infty_{n=1}
|a_n|^2\sum^\infty_{n=1}n^{-2\sigma},\qquad s=\sigma+it,$$
and so the abscissa of absolute convergence is at most $\frac12$ for the
series \thetag{2-1}. The abscissa of convergence may equal $\frac12$, as is
seen by the example $a_n=(n^{1/2}\log(n+1))^{-1}$.

As mentioned in the introduction, the collection of Dirichlet series
\thetag{2-1} satisfying \thetag{2-2} is denoted by $\calH$. It is a
complex Hilbert space when endowed with the inner product
$$\langle f,g\rangle_\calH=\sum^\infty_{n=1}a_n\bar{b}_n,$$
where $f(s)=\sum_n a_nn^{-s}$, and $g(s)=\sum_n b_nn^{-s}$. Thus, formally,
$\calH$ is just $l^2(\N)$, where $\N=\{1,2,3,\dots\}$ is the collection of
natural numbers. However, as a Hilbert space of analytic functions in the
half-plane $\Re s>\frac12$, $\calH$ has a rich and interesting structure. To
be specific, let us mention the classical problem of multiplication of two
Dirichlet series. Formally, the product of  $f(s)=\sum_n a_nn^{-s}$ and
$g(s)=\sum_n b_nn^{-s}$ is again a Dirichlet series
$$f(s)g(s)=\sum_n c_nn^{-s}, \tag 2-3$$
with
$$ c_n=\sum_{k,l:kl=n}a_kb_l. $$
For $f\in\calH$ and $g\in\calH$, the product $fg$ need not be in $\calH$, 
although the abscissa of absolute convergence of \thetag{2-3} can be shown 
to be at most $\frac12$ (as it is for functions in $\calH$). We are lead to 
the {\sl multiplier problem}: Find those functions $m$, analytic in $\Re s>
\frac12$, for which $m(s)f(s)$ is in $\calH$ for every $f\in\calH$. The 
collection of these {\sl multipliers} $m$ is denoted by $\calM$. Theorem 3.1
below solves the multiplier problem.

Let us finally mention a classical theorem which we will refer to from time
to time. We say that a set of real numbers $\xi_{1},\xi_{2},...,\xi_{r}$ is
{\sl $\Q$-linearly independent} if
$$c_{1}\xi_{1}+c_{2}\xi_{2}+\cdots+ c_{r}\xi_{r}=0$$
holds with integer coefficients $c_{1},c_{2},...,c_{r}$ only when the
coefficients are all zero. By the fundamental theorem of arithmetic, the prime
example of a $\Q$-linearly independent set of numbers is the image under the
logarithm function of any finite set of different primes. We have the
fundamental

\proclaim{Lemma 2.2} {\rm(Kronecker's theorem)}
Suppose the real numbers $1,\vartheta_{1},\vartheta_{2},...,\vartheta_{k}$
are $\Q$-linearly independent. Let $\alpha_{1},\alpha_{2},\ldots,\alpha_{k}$
be arbitrary real numbers, and $N$ and $\epsilon$ be given positive. Then
there are integers
$$n>N,\quad q_{1},q_{2},...,q_{k}$$
such that
$$|n\vartheta_{m}-q_{m}-\alpha_{m}|<\epsilon \quad (m=1,2,...,k).$$
\endproclaim

\noindent For a proof, we refer to Chapter XXIII of \cite{\bf\HW}.

To give an example of how Kronecker's theorem applies to our context, we
mention that an immediate consequence is the following identity:
$$\Big\|1+\sum_{p}a_{p}p^{-s}\Big\|_{\infty}=1+\sum_{p}|a_{p}|;$$
here the sum is only over prime indices, and
$$\|f\|_{\infty}=\sup_{\Re s>0}|f(s)|.\tag 2-4$$

\subhead 2.2. The Hardy space $H^{2}$ of the infinite-dimensional polydisk
\endsubhead Let $f$ be the Dirichlet series \thetag{2-1}, and suppose $f\in
\calH$. A fundamental observation, first made by Bohr \cite{\bf\Boo}, is that
if we put
$$z_1=2^{-s},z_2=3^{-s},\ldots\,,z_m=p_m^{-s},\ldots\,,$$
where $p_m$ denotes the $m$-th prime, then, in view of the fundamental
theorem of arithmetic, the Dirichlet series \thetag{2-1} can be considered
as a power series of infinitely many variables. 
The point is, as Bohr clarifies in his work, that the quantities 
$z_m=p^{-s}_m$, $m=1,2,3,\ldots$,
``obwohl sie ja in Wirklichkeit Funktionen nur des einen Parameters $s$ sind,
sich in mancher Beziehung fast ganz benehmen, als w\"{a}ren sie von einander
unabh\"{a}ngige Variable."
We put $z=(z_m)_m=(z_1,z_2,z_3,\ldots)$, and write
$$n=p^{\nu_1}_{k_1}p^{\nu_2}_{k_2}\ldots p^{\nu_r}_{k_r},$$
so that \thetag{2-1} formally takes the form
$$\fQ f(z)=\sum^\infty_{n=1}a_nz^{\nu_1}_{k_1}\cdot z^{\nu_2}_{k_2}\cdots
z^{\nu_r}_{k_r}.\tag 2-5$$
From now on, for a given element $f\in\calH$, $\fQ f$ denotes the corresponding
power series, and we drop the relationship between $z$ and $s$. 

Another, sometimes more convenient, way to think of the extension $\fQ f$ of
$f$ in \thetag{2-5} is to write it as
$$\fQ f(\phi)=\sum^\infty_{n=1}a_n\phi(n),\tag 2-5'$$
where $\phi$ is a {\sl quasi-character}, by which we mean a mapping $\phi:
\N\to\C$, which is multiplicative, $\phi(mn)=\phi(m)\phi(n)$, and has
$\phi(1)=1$ and $\phi(n)\in\D$ for $n>1$. Here, $\D$ is the open unit disk.
If we put, for the $m$-th prime $p_m$, $\phi(p_m)=z_m$, and identify $\phi$
with $z=(z_m)_m$, then \thetag{2-5} and \thetag{2-5'} express the same
function $\fQ f$. The Cauchy-Schwarz inequality applied to \thetag{2-5} (or
\thetag{2-5'}) and Euler's product formula give
$$\multline
|\fQ f(z)|^2=|\fQ f(\phi)|^2\leq\sum_{n=1}^{\infty}
|a_n|^2\sum^\infty_{n=1}|\phi(n)|^2\\
=\|f\|^2_{\calH}\cdot\prod_p(1-|\phi(p)|^2)^{-1}=\|f\|^2_{\calH}\cdot
\prod_{m=1}^\infty(1-|z_m|^2)^{-1},
\endmultline$$
provided that $|z_{m}|<1$ for every $m$. It follows that we have
bounded point evaluation at $z$, {\sl id est},
$$|\fQ f(z)|\leq C(z)\,\|f\|_{\calH},$$
if and only if $|z_m|<1$ for every $m$ and
$$\sum^\infty_{m=1}|z_m|^2<+\infty.$$

We denote by $\D^\infty$ the set of $z=(z_{m})_{m=1}^\infty$ for which
$|z_m|<1$ for every $m$, and call it {\sl the infinite-dimensional polydisk}.
The collection of power series \thetag{2-5} with $f\in\calH$ can be regarded
as a space of analytic functions on $\D^\infty\cap l^2({\Bbb N})$. We
denote this space by $H^2(\D^\infty)$, and supply it with the Hilbert space
structure induced by $\calH$. Note that $H^2(\D^\infty)$ differs from $\calH$ 
only in that its elements are considered as functions on $\D^\infty\cap 
l^2(\N)$ rather than on the half-plane $\Re s>\frac12$. In other words, 
$\calH$ results from considering \thetag{2-5'} for the special 
quasi-characters $\phi_s(n)=n^{-s}$.

Before discussing further the nature of $H^2(\D^\infty)$, we mention
an interesting problem, treated already by Hilbert \cite{\bf\Hi}: Is it
possible to extend the meaning of \thetag{2-5} beyond the set $\D^
\infty\cap l^2({\Bbb N})$, at least under certain favorable cirumstances?
As indicated by Bohr, Hilbert's solution to this problem has a remarkable
significance for Dirichlet series \cite{\bf\Boo}. Let
$$\Omega=\big\{z=(z_{m})_m\in\D^\infty:\,|z_m|<\beta_m\text{ for every }
m\big\}, $$
where $\{\beta_m\}_m$ is a sequence of positive numbers, $0<\beta_m\leq1$. For
a given point $z=(z_1,z_2,\ldots)\in\D^\infty$, we put
$$z^{(m)}=(z_1,z_2,\ldots,z_m,0,0,\ldots),$$
{\sl id est}, the $j$-th coordinate for $j>m$ is put equal to zero. For the
time being, let $f$ be a Dirichlet series \thetag{2-1}, about which we assume
only that it converges on some half-plane $\Re s>\sigma_c$ (we do not require
\thetag{2-2} to hold). The formal power series $\fQ f(z^{(m)})$ is the
``$m$-te Abschnitt" of $\fQ f$ at $z$, and it only depends on the finitely 
many coordinates $z_1,\ldots,z_m$. If the ``$m$-te Abschnitt" converges 
absolutely on $\Omega$ for each $m$, and there exists a constant $C=C(f)$ 
such that
$$|\fQ f(z^{(m)})|\leq C,\qquad z\in\Omega,$$
independently of $m$, $\fQ f$ is said to be {\sl finitely bounded} in $\Omega$.
In particular, the collection of series $\fQ f$ that are finitely bounded in
$\D^\infty$ is denoted by $H^\infty(\D^\infty)$, and we write
$$\|\fQ f\|_{H^\infty(\D^\infty)}=\sup_{m,z}|\fQ f(z^{(m)})|,$$
where $z$ ranges over $\D^\infty$. It is natural to define the value of
$\fQ f(z)$ as the limit of $\fQ f(z^{(m)})$ as $m\to+\infty$. The limit makes 
sense for $z\in\D^\infty\cap c_0(\N)$, where $c_0(\N)$ is the closed 
subspace of $l^\infty(\N)$ of sequences that converge to $0$. Moreover, the
function $\fQ f(z)$ is a bounded analytic function on $\D^\infty\cap c_0(\N)$.
Indeed, by the Schwarz lemma, we have, for $n<m$,
$$\big|\fQ f(z^{(m)})-\fQ f(z^{(n)})\big|\le 2\,\|\fQ f\|_{H^\infty(\D^\infty)}
\max\big\{|z_j|:\,n<j\le m\big\}.$$
For $z\in\D^\infty\cap c_0(\N)$, the maximum on the right hand side tends to
zero as both $n$ and $m$ tend to infinity, so that the above estimate says
that the $\fQ f(z^{(m)})$ form a Cauchy sequence in the space of bounded
analytic functions on $\D^\infty\cap c_0(\N)$. We note that the supremum of
$|\fQ f(z)|$ on $\D^\infty\cap c_0(\N)$ agrees with the norm of $f$ in
$H^\infty(\D^\infty)$.

It follows from our proof of the multiplier theorem in the next section that
$H^\infty(\D^\infty)$ corresponds to the set of Dirichlet series representing
bounded analytic functions in $\Re s>0$; in fact
$$\|f\|_\infty=\|\fQ f\|_{H^\infty(\D^\infty)}.\tag 2-6 $$
This norm identity reflects in a  concise way  
Bohr's observation that the $p^{-s}$ behave as ``independent variables''.

Thus far, we have thought about the space $H^{2}(\D^{\infty})$ as the Hardy
space $H^{2}$ on the infinite-dimensional polydisk $\D^{\infty}$. A perhaps
more natural setting is to regard the spaces $H^2(\D^\infty)$ and $H^\infty
(\D^\infty)$ as function spaces on the distinguished boundary $\T^\infty$. To
this end, we need a group-theoretical identification of $\T^\infty$, which we
shall now describe. Suppose that $\chi:{\Bbb N}\rightarrow{\Bbb C}$ satisfies
\roster
        \item"{(i)}"
        $\chi(mn)=\chi(m)\chi(n),\quad m,n=1,2,3,\ldots$
        \item"{(ii)}"
        $|\chi(n)|=1.$
\endroster
We say that $\chi$ is a {\sl character}, and write $\chi\in\Char$. We tacitly
understand that $\chi(n^{-1})=\chi(n)^{-1}$, so that the multiplicative formula
holds for all positive rational numbers $\Q_+$. {\sl The characters $\chi$
constitute the dual group of} $\langle\Q_+,\cdot\rangle$ ($\Q_+$ is given the
discrete topology, so that the space $\Char$ of characters is compact).
Examples of characters are the unit character $\chi(n)\equiv1$, and, more
generally, for $t\in{\Bbb R}$,
$$\chi(n)=n^{it}=e^{it\log n}.$$
The space $\Char$ can be identified with the infinite-dimensional polycircle
$\T^\infty$ in the following way. Given a point $z=(z_1,z_2,\ldots)\in
\T^\infty$, we define the value of $\chi$ at the primes through
$$\chi(2)=z_1,\quad\chi(3)=z_2,\ldots,\quad\chi(p_{m})=z_m,\ldots,$$
and extend the definition multiplicatively. This then yields a character, and
clearly, all characters are obtained by this procedure. This provides a
natural identification of the character space $\Char$ with $\T^\infty$.
The product topology on $\T^\infty$ makes it a compact space, and it is easily
shown that it corresponds to the topology $\Char$ gets as a dual group, so
that the identification $\Char=\T^\infty$ is topological.
There exists a unique (apart from normalization) Haar measure on $\Char$,
which we identify with the ordinary product measure $\rho$ on $\T^\infty$.
More precisely, let $\lambda$ be the arc length measure on $\T$, normalized so
that $\lambda(\T)=1$, and put, for Borel subsets $E_1,\,E_2,\ldots,\,E_N$ of
$\T$, and $E=E_1\times E_2\times\ldots\times E_N\times\T\times\T\times\ldots
\subset\T^\infty$,
$$\rho(E)=\lambda(E_1)\lambda(E_2)\cdot\ldots\cdot\lambda(E_N);$$
this defines a Borel measure $\rho$ on $\T^\infty$, which coincides with the
Haar measure on $\Char$, once it is agreed that the Haar measure should have
total mass $1$. We shall think of $\rho$ as living on $\Char$ as well as on
$\T^\infty$.

We return to the setting of the function $f\in\calH$ given by \thetag{2-1}.
The series
$$\sum_{n=1}^\infty a_{n}\chi(n)\tag 2-7$$
converges in the norm of $L^2(\Char,\rho)$ (though we do not know if we have
pointwise convergence almost everywhere) to a function $\P f(\chi)$, which is
determined uniquely by the two requirements that it be in $L^2(\Char,\rho)$
and that
$$\int_\Char\bar\chi(q)\P f(\chi)\,d\rho(\chi)=a_q,\qquad q\in\Q_+,$$
where we declare $a_q=0$ for $q\in\Q_+\setminus\N$. This function $\P f$ is an
extension to the characters of the earlier defined function $\fQ f$ on the
quasi-characters. Our shorthand for the above relationship between the
function $\P f(\chi)$ and the coefficients $\{a_n\}_n$ will be
$$\P f(\chi)\sim\sum_{n}a_{n}\chi(n),\qquad\chi\in\Char.$$
It is clear that the Plancherel identity 
$$\int_{\Char}|\P f|^2d\rho=\sum^\infty_{n=1}|a_n|^2=\|f\|^2_{\calH}$$
holds.

Let $H^2(\Char)$ (also written as $H^2(\T^\infty)$) be the closed subspace of
$L^2(\Char,\rho)$ consisting of functions that can be expanded in a series
\thetag{2-7}, with $\{a_n\}_n\in l^2(\N)$. Furthermore, let $H^\infty(\Char)$
be the intersection $L^\infty(\Char,\rho)\cap H^2(\Char)$, which is a closed
subspace of $L^\infty(\Char,\rho)$. To get the connection between the Hardy
spaces on the distinguished boundary $\T^\infty$ and on the interior $\D^
\infty$, we need the operator $\E:H^2(\T^\infty)\to H^2(\D^\infty)$, which
associates with a function $g(\chi)\sim\sum_n b_n\chi(n)$ in $H^2(\Char)$ the
function $\E g(\phi)=\sum_n b_n\phi(n)$ for quasi-characters $\phi\in\D^\infty
\cap l^2(\N)$. It is clear that it is an isometric isomorphism. The operators
$\fQ$, $\P$, and $\E$ are related via $\fQ=\E\P$.

\proclaim{Lemma 2.3} The operators $\P:\calH\to H^2(\Char)$ and $\E:
H^2(\Char)\to H^2(\D^\infty)$ are isometric isomorphisms. Furthermore, the
restriction of $\E$ to $H^\infty(\Char)$ is an isometric isomorphism
$H^\infty(\Char)\to H^\infty(\D^\infty)$.
\endproclaim

\demo{Proof} The first part just restates what was already done above. For the
rest, we can refer to \cite{\bf\CG}, where it is shown that the space
$H^\infty(\rho)$ is canonically isometrically isomorphic to both $H^\infty$
on $\D^\infty\cap c_0(\N)$ and on $\D^\infty\cap l^2(\N)$. In \cite{\bf\CG},
$H^\infty(\rho)$ is defined as the weak-star closure in $L^\infty(\T^\infty,
\rho)$ of the infinite-dimensional polydisk algebra $A(\D^\infty)$. We have
the inclusion $H^\infty(\rho)\subset H^2(\T^\infty)$, and by inspection, the
above canonical mapping coincides with our $\E$. It follows that it
suffices to check that the spaces $H^\infty(\rho)$ and $H^\infty(\T^\infty)$
are the same as (closed) subspaces of $L^\infty(\Char,\rho)$. Since both
$H^\infty(\rho)$ and $H^\infty(\T^\infty)$ are subspaces of $H^2(\T^\infty)$,
it is enough to check that their images under $\E$ coincide. In
\cite{\bf\CG}, it was shown that the operator $\E$ was in fact furnished by
integrating against the Poisson kernel in each variable, so that $\E$
applied to $H^\infty(\T^\infty)$ consists of {\sl bounded} analytic functions
on $\D^\infty\cap l^2(\N)$. By the above-mentioned result from \cite{\bf\CG},
the assertion of the lemma follows.
\qed\enddemo

\subhead 2.3. Vertical limit functions\endsubhead The notion of characters
allows us to clarify an important property of $\calH$. A natural unitary
operator on $\calH$ is that of vertical translation,
$$T_{\tau}f(s)=f(s+i\tau).$$
Fix an $f\in\calH$. To every sequence $\tau_{n}$ of translations there exists
a subsequence, say $\tau_{n(k)}$, such that $T_{\tau_{n(k)}}f(s)$ converges
uniformly on compact subsets of the domain $\Re s>\frac12$ to a limit function,
say $\tilde{f}(s)$. We will say that $\tilde{f}$ is {\sl a vertical limit
function of} $f$. We have the following result.

\proclaim{Lemma 2.4} The vertical limit functions of the function $f\in\calH$
given by \thetag{2-1} coincide with the functions of the form
$$f_{\chi}(s)=\sum_{n=1}^{\infty}a_{n} \chi(n) n^{-s},\tag 2-8$$
$\chi$ being a character.
\endproclaim

\demo{Proof} It is clear that every vertical limit function is of the form
(2-8). The fact that every function of this form is a vertical limit function
is a consequence of Kronecker's theorem (Lemma 2.2).
\qed\enddemo

\head
3. The theorem on multipliers
\endhead

We put
$${\Bbb C}_+=\{s=\sigma+it\in \C:\ \sigma>0\},$$
and let, as usual, $H^{\infty}({\Bbb C_{+}})$ denote the set of bounded
analytic functions on ${\Bbb C}_{+}$. By ${\Cal D}$ we denote the collection
of functions $f$ which can be represented by a convergent Dirichlet series
$$ f(s)=\sum^\infty_{n=1}a_nn^{-s} $$
in some half-plane $\sigma_{c}<\Re s<+\infty$. A {\sl multiplier} $m$ on
$\calH$ is by definition a holomorphic function on the half-plane $\Re s>
\frac12$ with the property that $mf\in\calH$ whenever $f\in\calH$. By standard
functional analysis, the operation of multiplying by a multiplier is a bounded
operator on $\calH$, and the multipliers form a commutative Banach algebra.
The collection of all multipliers on $\calH$ is denoted by $\calM$. We find it
convenient to consider a multiplier both as a function and as a bounded linear
operator on $\calH$. The operator norm of $m$ is denoted by $\|m\|_{\calM}$
and the supremum norm in $\C_{+}$ by $\|m\|_{\infty}$.

Our theorem on multipliers can be stated as follows.

\proclaim{Theorem 3.1} We have $\calM={\Cal D}\cap H^\infty(\C_+)$. Moreover,
$\|m\|_{\calM}=\|m\|_\infty$ holds for $m\in\calM$.
\endproclaim

The proof of Theorem 3.1 splits into two parts.

\subhead 3.1. Proof that $\calM\subset{\Cal D}\cap H^\infty(\C_+)$ and
$\|m\|_\infty\leq\|m\|_{\calM}$ \endsubhead
Let $m\in\calM$ be given. The idea of the proof is to show that $m$ lifts to a
multiplier on $H^2(\D^\infty)$, from which we deduce that it belongs to
$H^\infty(\D^\infty)$. This then entails that $m$ is a Dirichlet series, and
that it is bounded throughout $\C_+$. Since $1\in\calH$, we have $m\in\calH$,
and hence
$$m(s)=\sum^\infty_{n=1}b_nn^{-s},\qquad\Re s>{\tsize\frac12},$$
where the sequence $\{b_n\}_n$ is in $l^2(\N)$. So, if we apply the operator
$\P$ to $m$ and $f$, where $f\in\calH$, we get that both $\P m$ and $\P f$ are
in $H^2(\Char)$. Consequently, their product $\P m\,\P f$ is in $L^1(\Char,
\rho)$. Since $m$ is a multiplier, we have $mf\in\calH$, so that $\P(mf)\in
H^2(\Char)$. We wish to prove that $\P(mf)=\P m\,\P f$ as functions in $L^1(
\Char,\rho)$. To this end, note that for {\sl finite} Dirichlet series $f$,
this is verified by direct calculation. The general case when $f\in\calH$ is
arbitrary then follows by approximating $f$ with finite Dirichlet series.

Since $\|mf\|_{\calH}\le\|m\|_{\calM}\|f\|_{\calH}$, we get, by successively
plugging in $f=1,m,m^2$, and so on, that $\|m^j\|_{\calH}\le\|m\|_{\calM}^j$ holds
for $j=1,2,\ldots$. By what we just did, $\P(m^j)=(\P m)^j$, so that
$$\left(\int_{\Char}|\P m|^{2j}d\rho\right)^{1/(2j)}=\|(\P m)^j\|_{H^2(\Char)}
^{1/j}=\|m^j\|_{\calH}^{1/j}\le\|m\|_{\calM}.$$
As $j\to+\infty$, the left hand side tends to $\|\P m\|_{L^\infty(\Char,
\rho)}$. We conclude that $\P m$ belongs to $L^\infty(\Char,\rho)$ as well as
$H^2(\Char)$, that is, $\P m\in H^\infty(\Char)$. It also follows that
$\|\P m\|_{H^\infty(\Char)}\le \|m\|_{\calM}$. By Lemma 2.3, the function $\fQ
m=\E\P m$ is in $H^\infty(\D^\infty)$, and its norm equals the norm of $\P m$
in $H^\infty(\Char)$. Recall that by Hilbert's approach to power series of
infinitely many variables, the functions in the space $H^\infty(\D^\infty)$
are bounded and analytic in $\D^\infty\cap c_0(\N)$, and that the formula
giving $\fQ m$ on (part of) the quasi-characters is
$$\fQ m(\phi)=\sum_{n=1}^\infty b_n\phi(n),\qquad\qquad \phi\in\D^\infty\cap
l^2(\N).$$
Plugging in the special quasi-character $\phi_s(n)=n^{-s}$, we get our
function $m(s)$ back: $m(s)=\fQ m(\phi_s)$ for $\Re s>\frac12$. Since $\phi_s$
depends analytically on the parameter $s$ and is in $\D^\infty\cap c_0(\N)$
for $\Re s>0$, the fact that $\fQ m$ is bounded and holomorphic on $\D^\infty
\cap c_0(\N)$ implies that $m(s)$ is bounded and analytic in $\C_+$, and
$\|m\|_\infty\le\|\fQ m\|_{H^\infty(\D^\infty)}=\|\P m\|_{H^\infty(\Char)}
\le\|m\|_{\calM}$.

\subhead 3.2. Proof that ${\Cal D}\cap H^\infty({\Bbb C}_+)\subset{\calM}$
and $\|m\|_{\calM}\leq\|m\|_\infty$ \endsubhead
The key to the proof of the converse relation is the following lemma, which is
due to Fritz Carlson \cite{\bf\Car}.

\proclaim{Lemma 3.2} {\rm (Carlson's theorem)} Let $f(s)=\sum^\infty_{n=1}
a_nn^{-s}$ be convergent (and hence analytic) in $\C_+$ and bounded in every
half-plane $\Re(s)>\delta$ with $\delta>0$. Then, for  each $\sigma>0$,
$$\sum^\infty_{n=1}|a_n|^2n^{-2\sigma}=\lim_{T\to+\infty}\frac1{2T}
\int^T_{-T}|f(\sigma+it)|^2dt.$$
\endproclaim

An immediate corollary is the following.

\proclaim{Lemma 3.3} If $f(s)=\sum^\infty_{n=1}a_nn^{-s}$ is convergent and
bounded in $\C_+$, then $f\in\calH$ and
$$\|f\|_{\calH}=\lim_{\sigma\rightarrow 0^{+}}
\left(\lim_{T\to+\infty}\frac1{2T}\int^T_{-T}|f(\sigma+it)|^2dt
\right)^{1/2}.$$
\endproclaim

\demo{Remark} Note that together with Lemma 2.1, Lemma 3.3 implies Bohr's
inequality $\sigma_{a}-\sigma_{u}\leq\frac12$.
\enddemo

We carry on with the proof, and suppose that
$$m(s)=\sum^\infty_{n=1}b_nn^{-s}$$
converges in some half-plane $\Re s>\sigma_0$, and that it extends boundedly
and holomorphically to $\C_+$. Then, by Bohr's theorem (Lemma 2.1), the
Dirichlet series defining $m(s)$ actually converges uniformly to $m(s)$ on
every half-plane $\Re s>\vare$, with $\vare>0$. Let $f\in\calH$ have Dirichlet
series $\sum_{n=1}^\infty a_n n^{-s}$, and introduce, for $N=1,2,3,\ldots$,
the cut-off series
$$f_N(s)=\sum_{n=1}^N a_n n^{-s},\qquad s\in\C.$$
The function $mf_N$ is given by a convergent Dirichlet series in $\C_+$, and
 it is bounded there. We can now apply Lemma 3.3 to the function $mf_{N}$,
to obtain
$$\|mf_N\|_{\calH}\leq\|m\|_\infty\|f_N\|_{\calH}\leq\|m\|_\infty
\|f\|_{\calH}.$$
Since $\|mf\|_{\calH}\leq\sup_{N}\|mf_N\|_{\calH},$ it follows that
$$\|mf\|_{\calH}\leq\|m\|_\infty\|f\|_{\calH},$$
which completes the proof of Theorem 3.1.

\demo{Remark} In the proof of Theorem 3.1, we actually prove that the
multipliers on $\calH$ may be identified with the space $H^\infty(\D^\infty)$.
It is remarkable that just being able to extend the
function $m\in\calH$ holomorphically and boundedly to $\C_+$ should entail
that the function $\fQ m\in H^2(\D^\infty)$ is bounded on $\D^\infty\cap
l^2(\N)$ (and hence bounded and analytic on $\D^\infty\cap c_0(\N)$). After
all, the condition on $m$ just corresponds to the behavior of $\fQ m$ along the
one-dimensional complex variety which is the image of $\Cplus$ under $s\mapsto
\phi_s$.
\enddemo

\head
4. Some function theoretic properties of ${\Cal H}$
\endhead

\def\R{\Bbb R}

We first recall some results from ergodic theory. We then turn to the almost
sure behavior of the vertical limit functions $f_\chi$, with $\chi\in\Char$,
of a given $f\in\calH$. As an application, we consider the famous zeta
function, which here plays the role of the kernel function. After that, we
look at the function-theoretic properties of individual functions in $\calH$;
in particular, we study zero sets.

\subhead 4.1. Preliminaries from ergodic theory: Kronecker flows
\endsubhead Given a collection of real numbers $\alpha_1,\,\alpha_2,\alpha_3,
\ldots$, we consider the continuous group of vertical translations
$$T_\t(z_1,z_2,z_3,\ldots)=(e^{-i\t\alpha_1}z_1,e^{-i\t\alpha_2}z_2,
e^{-i\t\alpha_3}z_3,\ldots),$$
acting on the infinite-dimensional polycircle $\T^\infty$, where $\t$
ranges over the reals. The {\sl Kronecker flow} $\{T_\t\}_{\t}$ is known
to be ergodic if and only if for each fixed $n$, $n=1,2,3,\ldots$, the numbers
$\alpha_1,\,\alpha_2,\ldots,\,\alpha_n$ are $\Q$-linearly independent. This is
done for finite-dimensional polycircles $\T^n$ in \cite{{\bf\S}, pp. 64, 67,
69, 99}. After a few minor modifications, the proof in \cite{{\bf\S}}, which
is based on the approximation property in Kronecker's theorem (Lemma 2.2),
covers the infinite-dimensional case as well. We pick $\alpha_j=\log p_j$,
where $p_j$ denotes the $j$-th prime, and note that by the fundamental theorem
of arithmetic, any finite subset of the collection $\{\alpha_j\}_j$ is
$\Q$-linearly independent, so that the flow $\{T_\t\}_\t$ is ergodic. If we
write out the flow explicitly, we get
$$T_\t(z_1,z_2,z_3,\ldots)=\big(2^{-i\t}z_1,3^{-i\t}z_2,5^{-i\t}z_3,\ldots
\big),\qquad\t\in\R.$$
As in Section 2, points $z=(z_1,z_2,\ldots)$ in $\T^\infty$ are identified
with elements $\chi$ of the character group $\Char$ by putting $\chi(p_j)=
z_j$ for the $j$-th prime $p_j$. The flow then takes the more elegant form
$$(T_\t\chi)(n)=n^{-i\t}\chi(n),\qquad n\in\N,\,\,t\in\R.$$
We note that the ergodicity of the flow $\{T_t\}_t$ may now be checked off
directly from condition \thetag{iv} in \cite{{\bf\S}, p. 99}. By the
Birkhoff-Khinchin ergodic theorem \cite{{\bf\S}, pp. 11--12, 39, 99}, we have
$$\lim_{T\to+\infty}\frac1{2T}\int^T_{-T}g(T_\t\chi_0)\,d\t=
\int_{\Char}g(\chi)\,d\rho(\chi)\tag 4-1$$
for every $\chi_0\in\Char$ if $g$ is continuous on $\Char$, and for almost
every $\chi_0$ if we only assume that $g\in L^1(\Char,\rho)$. We now apply
this result in the context of our space $\calH$. For $f\in\calH$, with series
expansion $f(s)=\sum_n a_n n^{-s}$, we write $f_\sigma(s)=f(\sigma+s)=\sum_n
a_n n^{-\sigma-s}$, and note that for $\sigma\ge0$, this is again an element
of $\calH$. Recall that $\sigma_b(f)\in[-\infty,\frac12]$ is the abscissa of
boundedness for $f$, and that by Bohr's theorem (Lemma 2.1), it coincides with
the abscissa of uniform convergence. It follows that for $\sigma$,
$\sigma_b(f)<\sigma<+\infty$, the Dirichlet series for $f_\sigma(s)$ is
uniformly convergent on $\Re s>0$, so that by Kronecker's theorem (Lemma 2.2),
the partial sums of the series $\P f_\sigma$ converge uniformly on $\Char$,
making $\P f_\sigma$ continuous on $\Char$. Now, by \thetag{4-1},
$$\lim_{T\to+\infty}\frac1{2T}\int^T_{-T}|\P f_\sigma(T_\t\chi_0)|^2\,d\t=
\int_{\Char}|\P f_\sigma(\chi)|^2\,d\rho(\chi)=\sum_{n=1}^\infty|a_n|^2
n^{-2\sigma}\tag 4-2$$
holds for almost all $\chi_0\in\Char$ if $0\le\sigma\le\sigma_b(f)$, and for
all $\chi_0$ (in particular for the unit character $\chi_0\equiv1$) if
$\sigma>\frac12$. Notice the close resemblance with Carlson's theorem (Lemma
3.2). Indeed, Carlson's theorem can be read off from \thetag{4-2}, with $\chi_0
\equiv1$.

\subhead 4.2. The almost sure behavior of vertical limit functions\endsubhead
As before, let $f\in\calH$ be given by \thetag{2-1}, and write $f_\chi$ for
the vertical limit function of $f$, given by \thetag{2-8}.
The function $f_\chi$ shares with $f$ the property that it
is holomorphic on $\Re s>\frac12$. Helson \cite{\bf\Heb} has shown that the 
function $f_\chi$ extends analytically to $\Re s>0$ and that its Dirichlet 
series converges there, for almost every $\chi\in\Char$. We wish to illuminate
his argument, and obtain additional properties of the vertical limit functions.
To make the statement as precise as possible, we need the space $H^2_
{\text{i}}(\C_+)$ of functions $f$ holomorphic in $\C_+$ which have $f\circ
\varphi\in H^2(\D)$; here $H^2(\D)$ is the usual Hardy space of the unit disk
$\D$, and $\varphi$ is the Cayley transform $\varphi(z)=(1-z)/(1+z)$. The 
subscript ``i'' stands for (conformal) invariance. Let $\lambda_{\text{i}}$ 
be the probability measure $d\lambda_{\text{i}}(t)=\pi^{-1}(1+t^2)^{-1}$, 
which is obtained as the image under the Cayley transform of the normalized 
arc length measure on the unit circle:
$$\int_{-\infty}^{+\infty}|f(it)|^2d\lambda_{\text{i}}(t)=\frac1{2\pi}
\int_{-\pi}^\pi\big|f\circ\varphi(e^{i\theta})\big|^2d\theta,\qquad 
f\in H^2_{\text{i}}(\Cplus).$$
The space $H^2(\D)$ is usually regarded as a subspace of $L^2(\T)$, and
likewise the space $H^2_{\text{i}}(\C_+)$ may be considered to be a subspace of
$L^2(i\R,\lambda_{\text{i}}^*)$, the space of functions $g$ on $i\R$ that have 
$g(it)$ in $L^2(\R,\lambda_{\text{i}})$.

Apart from the ergodic statement, the following theorem is due to Helson 
\cite{{\bf \Heaa}}. For the benefit of the reader, we adapt Helson's proof to 
the present setting.
  
\proclaim{Theorem 4.1} Let $f\in\calH$ be given, with series expansion
\thetag{2-1}, and let $\varPi$ be a countable collection of absolutely 
continuous Borel probability measures on the real line. 
For almost every character $\chi\in\Char$, the function
$$f_\chi(s)=\sum_{n=1}^\infty a_n\chi(n)n^{-s},\qquad\Re s>{\tsize\frac12},$$
extends analytically to an element of $H^2_{\text{i}}(\C_+)$, has
$$\int_{-\infty}^{+\infty}\big|f_\chi(it)\big|^2d\varpi(t)<+\infty,\qquad
\text{for all}\quad\varpi\in\varPi,$$
and enjoys
$$\frac1{2T}\int_{-T}^T\big|f_\chi(it)\big|^2dt\to\sum_{n=1}^\infty|a_n|^2
\qquad\text{as}\quad T\to+\infty.$$
\endproclaim

\demo{Proof} We first decide on how to define $f_\chi(it)$. It should
correspond to the possibly divergent sum $\sum_n a_n\chi(n)n^{-it}$. This sum,
however, makes sense as $\P f(T_t\chi)$, for almost every $\chi$, where
$\{T_t\}_t$ is the ergodic Kronecker flow of the previous subsection. We thus
put $f_\chi(it) =\P f(T_t\chi)$, and observe that by \thetag{4-2}, with
$\sigma=0$, this function is locally in $L^2$ along the imaginary axis, and
has the asserted property
$$\frac1{2T}\int_{-T}^T|f_\chi(it)|^2dt\to\sum_{n=1}^\infty|a_n|^2\qquad
\text{as}\quad T\to+\infty,\tag 4-4$$
almost surely in $\chi$. Moreover, by Fubini's theorem, for a Borel 
probability measure $\varpi$,
$$\multline
\int_{\Char}\int_{-\infty}^{+\infty}\big|f_\chi(it)\big|^2d\varpi(t)\,d\rho
(\chi)=\int_{-\infty}^{+\infty}\int_{\Char}\big|f_\chi(it)\big|^2d\rho(\chi)\,
d\varpi(t)\\
=\int_{-\infty}^{+\infty}\|f\|_{\calH}^2\,d\varpi(t)=\|f\|_{\calH}^2<+\infty.
\endmultline$$
In particular, the function $f_\chi(it)$, considered as a function of $t$, 
almost surely is square integrable on the real line with respect to $\varpi$.
Elementary measure theory shows that the same holds true simultaneously for all
measures $\varpi$ in $\varPi$, since the latter set is countable. By the same
token, $f_\chi(it)$ is almost surely in $L^2(\R,\lambda_{\text{i}})$.

A function $g$ in $L^2(\T)$ is in $H^2(\D)$ if and only if 
$$\int_\T z^n g(z)\,d\lambda(z)=0,\qquad n=1,2,3,\ldots,$$
where $\lambda$ is as before the normalized arc length measure on the unit 
circle. After an application of the Cayley transform, we have that a function 
$g\in L^2(i\R,\lambda_{\text{i}}^*)$ is in $H^2_{\text{i}}(\C_+)$ if and only 
if 
$$\int_{-\infty}^{+\infty}\left(\frac{1-it}{1+it}\right)^n g(it)\,d\lambda_
{\text{i}}(t)=0,\qquad n=1,2,3,\ldots.$$ 
Therefore, to see that $f_\chi$ is almost surely in $H^2_{\text{i}}(\C_+)$,
it suffices to check that
$$\int_{-\infty}^{+\infty}\left(\frac{1-it}{1+it}\right)^n f_\chi(it)\,
d\lambda_{\text{i}}(t)=0,\qquad n=1,2,3,\ldots,$$ 
for almost all $\chi$. However, to check that an $L^2(\Char,\rho)$ function 
vanishes almost everywhere, it is enough to show that all of its Fourier 
coefficients are $0$. We set $a_q=0$ for $q\in\Q_+\setminus\N$, and integrate 
the left hand side of the above expression against $\bar\chi(q)$ to get the 
Fourier coefficients:
$$
\thickmuskip=.48 \thickmuskip
\medmuskip=.48 \medmuskip
\multline
\int_\Char\bar\chi(q)\int_{-\infty}^{+\infty}\left(\frac{1-it}{1+it}
\right)^n f_\chi(it)\,d\lambda_{\text{i}}(t)\\
=\int_{-\infty}^{+\infty}
\left(\frac{1-it}{1+it}\right)^n\int_\Char\bar\chi(q)\,f_\chi(it)\,
d\rho(\chi)\,d\lambda_{\text{i}}(t)
=\int_{-\infty}^{+\infty}\left(\frac{1-it}{1+it}\right)^n a_qq^{-it}\,
d\lambda_{\text{i}}(t)=0,
\endmultline$$ 
where we have used Fubini's theorem, and that $a_q=0$ for $q<1$. This completes
the proof.
\qed
\enddemo

The ergodic reasoning behind \thetag{4-4} leads to an estimate of $f_{\chi}$, 
which does not seem to follow from Helson's work.

\proclaim{Theorem 4.2} Let $f$ and $f_\chi$ be as in Theorem 4.1, and write
$s=\sigma+it$. Then almost surely in $\chi\in\Char$,
$$\big|f_\chi(s)-a_1\big|\le C\,\left(\frac{1+|t|^{1/2}}{\sigma^{1/2}}\right),
\qquad s\in\C_+,$$
for some constant $C=C(f,\chi)$, $0<C<+\infty$. Moreover, almost surely in
$\chi$,
$$f_\chi(s)=a_1+o\left(\frac{|t|^{1/2}}{\sigma^{1/2}}\right)\qquad\text{as}
\quad|t|\to+\infty$$
holds uniformly in $\sigma>0$.
\endproclaim 

\demo{Proof} As in \thetag{4-4}, we have, by ergodic theory, that almost 
surely in $\chi$,
$$\frac1T\int_0^T\big(f_\chi(i\tau)-a_1\big)^2d\tau\to0\qquad\text{as}\quad
|T|\to+\infty,$$
the space average of the function (in $\chi$) being $0$. This entails
that the function
$$F_{\chi}(ix)=\int_0^x\big(f_\chi(i\tau)-a_1\big)^2d\tau$$
meets  
$$F_\chi(ix)=o(|x|),\qquad\text{as}\quad|x|\to+\infty,\tag 4-5$$ 
almost surely in $\chi$. By Theorem 4.1, the function $f_\chi$ is almost 
surely in $H^2_{\text{i}}(\C_+)$, so that the squared function $(f_\chi-a_1
)^2$ almost surely belongs to the analogously defined space $H^1_{\text{i}}
(\C_+)$. As such it is given by the Poisson formula
$$\big(f_\chi(\sigma+it)-a_1\big)^2=\frac\sigma\pi\int_{-\infty}^{+\infty}
\,\frac{\big(f_\chi(ix)-a_1\big)^2}{\sigma^2+(t-x)^2}\,dx.$$
Integrating by parts, using \thetag{4-5}, we obtain the representation 
formula
$$\big(f_\chi(\sigma+it)-a_1\big)^2=\frac{2\sigma}\pi\int_{-\infty}^{+\infty}
\frac{t-x}{\big(\sigma^2+(t-x)^2\big)^2}\,F_\chi(ix)\,dx,
$$
for $\sigma>0$ and $t\in\R$. After an application of the size control 
\thetag{4-5} to this integral, the desired estimates follow by taking square 
roots. 
\qed
\enddemo

For a certain class of functions in $\calH$, the estimate of Theorem 4.2 can
be improved considerably. Note that the conclusion is that the growth in the
imaginary direction is precisely what the Schnee-Landau theorem requires to
imply convergence of Dirichlet series \cite{\bf\Booo{\rm,} \Sch}.

\proclaim{Corollary 4.3} Let $f\in\calH$ be such that all powers $f^N$,
where $N$ is a positive integer, also are in $\calH$. Then, for each $\vare>0$,
we have, almost surely in $\chi\in\Char$,
$$\big|f_\chi(s)-a_1\big|\le C\,\left(\frac{1+|t|^\vare}{\sigma^\vare}\right),
\qquad s=\sigma+it\in\C_+,$$
for some constant $C=C(\chi,f,\vare)$, $0<C<+\infty$. 
\endproclaim

\demo{Proof} The estimate is more or less immediate from Theorem 4.2. 
\qed
\enddemo

The general question about almost sure convergence of the Dirichlet series of
vertical limit functions was treated by Helson in \cite{\bf\Heb}, in a 
somewhat more general context. He obtained the following basic result.  

\proclaim{Theorem 4.4} {\rm (Helson)} Let $f\in\calH$ be given, with series 
expansion \thetag{2-1}. For almost every character $\chi$, the Dirichlet series
$$ f_\chi(s)=\sum_n a_n \chi(n) n^{-s} $$
converges in the half-plane $\Re s>0$.
\endproclaim

According to a classical formula for the abscissa of convergence 
\cite{{\bf\HR}, pp. 6--8} (see also Section 2 of the present paper), Theorem 
4.4 is equivalent to the statement that, almost surely in $\chi$, the partial 
sum function
$$S_N(\chi)=\sum_{n:n\leq N} a_n\chi(n)\tag4-6$$  
has $S_N(\chi)=O(N^\vare)$ as $N\to+\infty$, almost surely in $\chi$, for all
$\vare>0$. So, in light of  Helson's theorem, we find the following special 
case interesting. 

\proclaim{Theorem 4.5} Let $f\in\calH$ be given, with series expansion 
\thetag{2-1}. Suppose $a_n=0$ for all composite numbers. For almost every 
character $\chi$, we then have $S_N(\chi)=O(1)$ as $N\rightarrow \infty$, 
where $S_N(\chi)$ is as in \thetag{4-6}.
\endproclaim

\demo{Proof} The non-composite numbers are the primes and unity ($1$). Without
loss of generality, we may suppose that $a_1=0$.
As the parameter $p$ runs through the primes, the $\chi(p)$, treated as 
functions of $\chi$, run through the distinct coordinate variables in the
polycircle ${\Bbb T}^\infty$, which is our standard realization of the 
character space $\Xi$. This has a clear interpretation: {\sl the random
variables
$$
\chi\mapsto a_p \chi(p), \ \ \ p=2,3,5,7,11,...
$$
are mutually independent}, have mean squares 0, and variances 
$\sigma_p^2=|a_p|^2$.

For a given positive integer $N$ and a positive real number $M$, let $E(N,M)$
be the set of all characters $\chi$ for which $|S_N(\chi)|<M$. By 
Kolmogorov's inequality \cite{{\bf \S}, p. 260}, we have
$$\rho\big(\cap_NE(N,M)\big)\geq 1-M^{-2}\sum_p \sigma_p^2=1-M^{-2}|a_p|^2.$$
By letting $M$ tend to $+\infty$, we see that $S_N(\chi)=O(1)$ as 
$N\rightarrow +\infty$ holds for almost all $\chi\in\Xi$. \qed
\enddemo

\subhead 4.3. The kernel function. Riemann's zeta function for random
characters\endsubhead
Gi\-ven a separable Hilbert space $\calA$ of analytic functions on a domain
$\Omega$ in $\C$ or in $\C^n$, one forms the kernel function $K_{\calA}(z,w)$ 
by taking some orthonormal basis $\{e_n(z)\}_n$ in $\calA$, and putting
$$K_{\calA}(z,w)=\sum_n e_n(z)\,\bar e_n(w),\qquad (z,w)\in\Omega\times
\Omega.$$
This then proves to be independent of the particular choice of orthonormal
basis, and has the reproducing property that
$$f(w)=\langle f,K_{\calA}(\cdot,w)\rangle_{\calA},\qquad w\in\Omega.$$
In fact, as an element of $\calA$, $K(\cdot,w)$ is uniquely determined by its
reproducing property. In $\calH$, an orthonormal basis is supplied by $e_n(z)
=n^{-z}$, for $n=1,2,3,\ldots$, so that its kernel function is
$$K_{\calH}(z,w)=\sum_{n=1}^\infty n^{-z-\bar w}=\zeta(z+\bar w),\qquad
\Re z>{\tsize\frac12},\,\,\Re w>{\tsize\frac12},$$
where $\zeta(s)$ is the Riemann zeta function:
$$\zeta(s)=\sum_{n=1}^\infty n^{-s},\qquad\Re s>1.$$
In Section 2, we modeled the space $\calH$ as both the Hardy space on the
infinite-dimensional polydisk $\D^\infty$ and the polycircle $\T^\infty$. The
polycircle $\T^\infty$ was identified with the character group $\Char$ of the
multiplicative positive rationals. In particular, for quasi-characters $\phi
\in\D^\infty\cap l^2(\N)$, the "point evaluation" $f\mapsto \fQ f(\phi)$ is
a continuous linear functional on $\calH$, so that it too must be given by
a kernel function $K_{H^2(\D^\infty)}(\psi,\phi)$,
$$\fQ f(\phi)=\langle\fQ f,K_{H^2(\D^\infty)}(\cdot,\phi)
\rangle_{H^2(\D^\infty)},\qquad\phi\in\D^\infty\cap l^2(\N).$$
In terms of a series expansion, it is written
$$K_{H^2(\D^\infty)}(\psi,\phi)=\sum_{n=1}^\infty\psi(n)\,\bar\phi(n),\qquad
\psi,\phi\in\D^\infty\cap l^2(\N).$$
The inner product in $\calH$, and hence in $H^2(\D^\infty)$, is better
visualized on the distinguished boundary $\T^\infty=\Char$, where we have
$$\multline
\langle f,g\rangle_{\calH}=\langle\fQ f,\fQ g\rangle_{H^2(\D^\infty)}
=\langle\P f,\P g\rangle_{H^2(\T^\infty)}\\
=\int_\Char\P f(\chi)\,\overline
{\P g}(\chi)\,d\rho(\chi),\qquad f,g\in\calH.
\endmultline$$
This suggests introducing the kernel
$$K_{H^2(\Char)}(\chi,\phi)=\sum_{n=1}^\infty\chi(n)\,\bar\phi(n),\qquad
\chi\in\Char,\,\,\phi\in\D^\infty\cap l^2(\N),$$
where the sum, for fixed $\phi$, is understood to converge in the sense of
the space $L^2(\Char,\rho)$. This kernel then has the reproducing property
$$\multline
\fQ f(\phi)=\langle\P f,K_{H^2(\Char)}(\cdot,\phi)
\rangle_{H^2(\Char)}\\
=\int_\Char \P f(\chi)\,\overline{K_{H^2(\Char)}}
(\chi,\phi)\,d\rho(\chi),\qquad\phi\in\D^\infty\cap l^2(\N).
\endmultline$$
On the other hand, it is well known that the reproducing kernel on $H^2(\T^
N)$ for finite-dimensional polydisks is given as the product of the Cauchy
kernel in each variable,
$$K_{H^2(\T^N)}(z,w)=\prod_{n=1}^N (1-\bar w_nz_n)^{-1},\qquad z\in\T^N,\,\,
w\in\D^N,$$
and in the limit as the dimension $N$ tends to infinity, we get \cite{\bf\CG}
$$\multline
K_{H^2(\T^\infty)}(\chi,\phi)=\prod_{n=1}^\infty(1-\bar w_nz_n)^{-1}\\
=\prod_p \big(1-\bar\phi(p)\chi(p)\big)^{-1},\qquad
\chi\in\T^\infty,\,\,\phi\in\D^\infty\cap l^2(\N),
\endmultline$$
where the second product runs over the primes, and $z_n=\chi(p_n)$, $w_n=
\phi(p_n)$, for the $n$-th prime $p_n$. For $\phi\in\D^\infty\cap l^1(\N)$,
the above product converges pointwise to a continuous function of $\chi$,
but in general we must interpret the above product as being convergent in
$L^2(\Char,\rho)$ \cite{\bf\CG} (Cole and Gamelin used martingale theory to
obtain the convergence). Since the kernel function of a given Hilbert space of
analytic functions is unique, we arrive at the equality $K_{H^2(\Char)}=K_
{H^2(\T^\infty)}$, that is, {\sl the Euler identity}
$$\sum_{n=1}^\infty\chi(n)\,\bar\phi(n)=\prod_p \big(1-\chi(p)\bar\phi(p)
\big)^{-1},\qquad\chi\in\T^\infty,\,\,\phi\in\D^\infty\cap l^2(\N),$$
which holds pointwise in $\chi$ for $\phi\in\D^\infty\cap l^1(\N)$, and almost
everywhere for general $\phi$. This is a surprising interpretation of the Euler
identity as arising from two ways of looking at one and the same kernel
function! If we specialize to the particular quasi-characters $\phi(n)=
n^{-\bar s}$, we get the more familiar
$$\sum_{n=1}^\infty\chi(n)\,n^{-s}=\prod_p \big(1-\chi(p)p^{-s}\big)^{-1},
\qquad\chi\in\Char,\,\,\Re s>{\tsize\frac12},\tag 4-7$$
with pointwise convergence in $\chi$ for $\Re s>1$. We shall write $\zeta_\chi
(s)$ for the analytic function on $\Re s>1$ given by either side of
\thetag{4-7}. By Lemma 2.4, suitably modified, these are the {\sl vertical
limit functions of the zeta function} $\zeta(s)$.

\proclaim{Theorem 4.6} Suppose the coefficients $\{a_n\}_n$ are totally
multiplicative and square summ\-able, with $a_1=1$. Then, for almost every
character $\chi$, the Dirichlet series
$$f_\chi(s)=\sum_{n=1}^\infty a_n\chi(n)n^{-s},\qquad\Re s>{\tsize\frac12},$$
converges to a zero-free analytic function on the half-plane $\Re s>0$.
\endproclaim

\demo{Proof} The convergence statement follows from Theorem 4.4. We check that
the function $f_\chi(s)$ almost surely lacks zeros in $\Re s>0$. Let $\mu(n)$ 
be the M\"obius function, which has $\mu(1)=1$, $\mu(n)=(-1)^k$ if $n$ is the 
product of $k$ different primes, and $\mu(n)=0$ if $n$ is divisible by a 
square (other than $1$). The M\"obius function enables us to express the 
reciprocal of $f_\chi(s)$,
$$1/f_{\chi}(s)=\sum_{n=1}^{\infty}\mu(n)a_{n}\chi(n)n^{-s},\qquad\Re s>1,$$
because the coefficients $\{a_n\}_n$ are totally multiplicative. By Theorem
4.1, $1/f_\chi(s)$ extends analytically to $\C_+$ almost surely in $\chi$, 
whence the assertion follows.
\qed
\enddemo

Theorem 4.6 has a curious interpretation for the vertical limit functions
$\zeta_\chi(s)$ of the Riemann zeta function: the Riemann hypothesis holds for
almost every $\chi$ (see also \cite{\bf\Heb}).

\proclaim{Corollary 4.7} {\rm (Helson)} For almost every character $\chi$, 
the Dirichlet series $\zeta_\chi(s)=\sum_n\chi(n)n^{-s}$ converges on the 
half-plane $\Re s>\frac12$ to an analytic function which has no zeros there.
\endproclaim

\demo{Proof} Apply Theorems 4.1 and 4.6 to the coefficients $a_n=n^{-\frac12
-\vare}$, with $\vare>0$. The assertion follows, but only in the slightly
smaller half-plane $\Re s>\frac12+\vare$. An elementary measure-theoretical 
argument now yields the desired result.
\qed
\enddemo

By the formula for the abscissa of convergence of a Dirichlet series (see 
Section 2), Corollary 4.7 has the following consequence.

\proclaim{Corollary 4.8} Let $\vare>0$. For almost every character $\chi$, we
have $\sum_{n=1}^N\chi(n)=O\big(N^{1/2+\vare}\big)$ as $N\to+\infty$.
\endproclaim

\demo{Question} Can the above be sharpened to the statement that, almost 
surely, \break
$\sum_{n=1}^N\chi(n)=O\big(\sqrt{N}\log N\big)$ as $N\to+\infty$?
\enddemo

Using Corollary 4.3, we can control the growth of the logarithm of $\zeta
_\chi(s)$. The branch of the logarithm intended is the one with
$$\log\zeta_\chi(s)=-\sum_p\log\big(1-\chi(p)p^{-s}\big),\qquad \Re s>1,$$
where the logarithms on the right hand side are given by the principal branch.

\proclaim{Theorem 4.9} Suppose the coefficients $\{a_n\}_n$ are totally
multiplicative and square summ\-able, with $a_1=1$. Then, for almost every
character $\chi$, the function $f_\chi(s)=\sum_{n}a_n\chi(n)n^{-s}$ has a 
logarithm
$$\log f_\chi(s)=-\sum_p\log\big(1-a_p\chi(p)p^{-s}\big),\qquad \Re s>{\tsize
\frac12},$$
(the right hand side involves the principal branch) which extends 
holomorphically to $\Re s>0$, and enjoys $(s=\sigma+it)$
$$\log f_\chi(s)=o\left(\log\frac{|t|}{\sigma}\right)\qquad\text{as}\quad
|t|\to+\infty,$$
uniformly in $\sigma$, $\sigma>0$.
\endproclaim 

\demo{Proof} Let $\alpha$ be a complex parameter. We define the power
$$f(s)^\alpha=\prod_p\big(1-a_p p^{-s}\big)^{-\alpha},\qquad\Re s>{\tsize\frac
12},$$
where the right hand side employs the principal branch of the logarithm.
The identification of $\calH$ with $H^2(\Char)$ in Section 2 allows us to 
express the square of the norm of $f^\alpha$ in $\calH$ as 
$$\prod_p\frac1{2\pi}\int_{-\pi}^\pi\big|\big(1-a_pe^{i\theta}\big)^{-\alpha}
\big|^2d\theta.$$
Using a Maclaurin expansion of $\big(1-a_pe^{i\theta}\big)^{-\alpha}$, one 
sees that for $\alpha$ confined to a compact subset of $\C$,
$$\frac1{2\pi}\int_{-\pi}^\pi\big|\big(1-a_pe^{i\theta}\big)^{-\alpha}\big|^2
d\theta=1+O(|a_p|^2).$$
It follows that $f^\alpha\in\calH$. In particular, $f^{\gamma N}\in\calH$ for
all $N=1,2,3,\ldots$, and $\gamma=1,-1,i,-i$. By Corollary 4.3, we have, for
each $\vare>0$, and almost all $\chi$,
$$\big(f_\chi(s)\big)^\gamma=O\left(\frac{|t|^\vare}{\sigma^\vare}\right)\qquad
\text{as}\quad|t|\to+\infty,$$
uniformly in $\sigma$, $\sigma>0$, whence the assertion follows.
\qed
\enddemo

\proclaim{Corollary 4.10} Fix $\vare>0$. Almost surely in $\chi$, 
$$\log\zeta_\chi(\sigma+it)=o(\log|t|)\qquad \text{as}\quad|t|\to+\infty$$
holds uniformly in $\sigma$, $\sigma>\frac12+\vare$.
\endproclaim

\demo{Proof} Apply Theorem 4.9 to the coefficients $a_n=n^{-\frac12-
\delta}$, with $\delta=\frac12\vare$. 
\qed
\enddemo

\subhead 4.5. Function-theoretic properties of individual functions in
$\calH$ \endsubhead
The purpose of this paragraph is to obtain some basic information about the
function theoretic properties of individual functions in $\calH$. In
particular, we are interested in the structure of their zero sequences. The
situation starkly contrasts the "almost sure" behavior that we have
concentrated on so far, because here, the functions are only known to be
analytic in the half-plane $\Re s>\frac12$ (which we denote by $\C_{\frac12}$),
and zeros may actually accumulate at a boundary point such as $s=\frac12$.

Let $H^2_\infty(\C_{\frac12})$ denote the uniformly local $H^2$ space on $\C_
{\frac12}$: a function $g$ holomorphic on $\C_{\frac12}$ is said to be in it
if
$$\sup_{\theta\in\R}\,\sup_{\sigma>1/2}\int_\theta^{\theta+1}|g(\sigma+it)
|^2dt<+\infty.$$
It is a Banach space, and the functions in it are bounded in every half-plane
$\Re s>\sigma_0$, with $\sigma_0>\frac12$.

\proclaim{Theorem 4.11} We have the inclusion $\calH\subset H^2_\infty
(\C_{\frac12})$, and the injection mapping is continuous.
\endproclaim

\demo{Proof} Let $f\in\calH$ have the series expansion \thetag{2-1}. We wish
to prove that for all $\sigma>\frac12$,
$$\int^{\theta+1}_\theta|f(\sigma+it)|^2dt\leq C\sum^\infty_{n=1}|a_n|^2\tag
4-8$$
holds, where $C$ is an absolute constant. We note that it suffices to obtain
\thetag{4-8} for finite Dirichlet series $f$, since on compact subsets of $\C_
{\frac12}$, elements of $\calH$ are uniformly approximable by them. Moreover,
by making a vertical translation, we see that we can set $\theta=0$, and by
the Poisson integral formula, we see that it suffices to consider the limit
case $\sigma=\frac12$. By duality, we have
$$\multline
\int^1_0|f({\tsize\frac12}+it)|^2dt=\sup_g\Big|\int_0^1
\sum_{n} a_n n^{-1/2}e^{-it\log n}g(t)\,dt\Big|^2\\
=\sup_g\Big|\sum_{n} a_n n^{-1/2}\widehat g(\log n)\Big|^2
\le\sum_{n=1}^\infty|a_n|^2\,\sup_g\sum_{n=1}^\infty n^{-1}|\widehat g(\log
n)|^2,
\endmultline\tag 4-9$$
where the supremum is taken over all $g\in L^2(0,1)$ of norm $\le1$, and
$\widehat g$ is the Fourier transform of $g$,
$$\widehat g(\xi)=\int_0^1e^{-it\xi}g(t)\,dt,$$
which extends to an entire function of exponential type $\le1$. For such
functions $\widehat g$ ,
$$\sum_{n=1}^\infty n^{-1}|\widehat g(\log n)|^2\le C\int_{-\infty}^{+\infty}
|\widehat g(\xi)|^2d\xi$$
holds for some absolute constant $C$. This follows from a suitable adaptation
of a classical inequality of Plancherel and P\'{o}lya \cite{{\bf\Y}, pp.
96--98}. Heuristically speaking, the reason is that the left hand side looks
like a Riemann sum of part of the integral to the right, and that the
functions $\widehat g$ are sufficiently smooth. Modulo the Plancherel identity,
the assertion now follows from \thetag{4-9}.
\qed
\enddemo

\proclaim{Corollary 4.12} Let $f\in\calH$, $f(s)\not\equiv0$. If $\{s_k\}_k$,
$s_k=\sigma_k+it_k$, denotes the sequence of zeros of $f$ in $\C_{\frac12}$,
then
$$\sup_{\theta\in\R}\sum_{k:s_k\in Q_\theta}(\sigma_k-{\tsize\frac12})
<+\infty,$$
where $Q_\theta$ is the half-strip $\Re s>0$, $\theta<\Im s<\theta+1$.
\endproclaim

\demo{Proof} By Theorem 4.11, $f(s)/s$ is contained in $H^2(\C_{\frac12})$, and
so $\{s_k\}_k$ satisfies the Blaschke condition
$$\sum_k\frac{\sigma_k-{\frac12}}{1+|s_k|^2}<+\infty.$$
All convergent nontrivial Dirichlet series have a zero-free half-plane, so that
$$\sup_k\sigma_k<+\infty.$$ Hence
$$A(\theta)=\sum_{k:s_k\in Q_\theta}(\sigma_k-{\tsize\frac12})<+\infty,$$
for each real $\theta$. If the function $A(\theta)$ were unbounded, we could
find a sequence $\theta_k$ of reals such that $A(\theta_k)\to+\infty$. But
then the functions $f_k(s)=f(s+i\theta_k)$ would tend to $0$ uniformly on
compact subsets of $\C_{\frac12}$, because they would have an ever increasing
mass of zeros, and their norms in $H^2_\infty(\C_{\frac12})$ are uniformly
bounded. By Lemma 2.4, any vertical limit function has the form
$$\sum^\infty_{n=1}a_n\chi(n)n^{-s},$$
where $\chi$ is a character. Such a function cannot be identically zero. This
is a con\-tra\-dic\-tion, and so $A(\theta)$ is bounded.
\qed
\enddemo

\proclaim{Proposition 4.13} There exists a function $f\in\calH$, other than the
$0$ function, whose zero set contains a subsequence tending to $s=\frac12$
along the real line.
\endproclaim

\demo{Proof} Let $\{b_n\}_n$ be the sequence $b_1=0$, $b_n=n^{-1/2}(\log
n)^{-1}$ for $n>1$, and notice that it has $\sum_n b_n^2<+\infty$ and $\sum_n
n^{-1/2}b_n=+\infty$. We now select integers $n_k$ and points $\sigma_k$
according to the following fashion. Let $n_0=1$. Given $n_k$, we pick
$\sigma_k>\frac12$ so close to $\frac12$ and $n_{k+1}$ so much larger than
$n_k$ that
$$\sum_{n=n_k}^{n_{k+1}-1}n^{-\sigma_k}b_n-\sum_{n=1}^{n_k-1}n^{-\sigma_k}b_n
-\sum_{n=n_{k+1}}^{\infty}n^{-\sigma_k}b_n>0.$$
%(For example, $n=2^{2^{2^{k}}}$ and $\sigma_{k}=\frac12+2^{-2^{k+1}}$
%will do.)
The sequence $\{\sigma_k\}_k$ is strictly decreasing, with limit $\frac12$.
If we put $a_n=(-1)^kb_n$ for $n_k\le n<n_{k+1}$, and let $f$ be as in
\thetag{2-1}, then $f$ is real-valued on the real half-axis $]\frac12,+
\infty[$, and the sign of $f(\sigma_k)$ alternates as $(-1)^k$.
By the intermediate value theorem of calculus, $f$ has infinitely many zeros
along the real line, which accumulate at $s=\frac12$.
\qed
\enddemo

% \head
% 6. Cyclicity in ${\Cal H}$ and ${\Cal M}$
% \endhead

\head
5. Systems of dilated functions in $L^2(0,1)$
\endhead

\define\var{\varphi}

In this section, we return to the basis and completeness problems described
in the introduction.

Let $\var\in L^{2}(0,1)$ be given and again consider it as an odd periodic
function of period 2.
The  basic questions are: For which $\var$ is the system
$$ \var(x),\;\var(2x),\;\var(3x),\ldots$$
a Riesz basis in $L^2(0,1)$, and for which $\var$ is the same system complete
in $L^2(0,1)$? We recall that completeness means that the system spans a dense
linear subspace. The observation that the only orthogonal bases of this kind
are generated by the functions $\var(x)=C\sin(\pi x)$ motivates introducing the
representation
$$\var(x)=\sum_{n=1}^{\infty}a_{n}\sqrt{2}\sin(n\pi x)$$
for an arbitrary $\var\in L^{2}(0,1)$, and seeking solutions in terms of the
coefficients $a_n$. For both problems, a necessary condition is $a_1\neq 0$,
because otherwise the function $\sin(\pi x)$ cannot be approximated. It is no
restriction to normalize: $a_1=1$ from now on.

Let $\{e_n\}_n$ be the standard orthonormal basis $e_n(x)=\sqrt{2}\sin(n\pi
x)$. Following Beurling \cite{\bf\Be}, we associate to each $f\in L^2(0,1)$
with sine series expansion $f(x)=\sum_n c_ne_n(x)$ the auxiliary Dirichlet
series
$$Sf(s)=\sum^\infty_{n=1}c_nn^{-s}.$$
Let $\|\cdot\|_{L^2}$ and $\langle\cdot,\cdot\rangle_{L^2}$ be the norm and
inner product in $L^2(0,1)$, respectively. The operator $S$ is an isometric
isomorphism between $L^2(0,1)$ and $\calH$: $\|Sf\|_{\calH}=\|f\|_{L^2}$. Also,
inner products are preserved: $\langle Sf,Sg\rangle_{\calH}=\langle f,g
\rangle_{L^2}$. Given a function $f$ in $L^2(0,1)$ of the form $f(x)=\sum_n
c_n e_n(x)$, where only finitely many $c_n$ are nonzero, we associate to it
the function $T_{\var}f(x)=\sum_n c_n \var_n(x)$, where $\var_n(x)=\var(nx)$.
It is clear that $T_{\var}f\in L^2(0,1)$, given the restriction on $f$. When
the operator $S$ is applied to $T_{\var}f$, we arrive at the identity
$$S(T_{\var}f)(s)=S\var(s)\,Sf(s),\qquad \Re s>{\tsize\frac12}.\tag{5-1}$$
This has the interpretation that replacing the basis sequence $\{e_n\}_n$ with
$\{\var_n\}_n$ corresponds to multiplication by $S\var$ on the $S$-transformed
side.

\subhead 5.1 The function $1/S\var$. The biorthogonal system\endsubhead
Let $S\var(s)=\sum_n a_nn^{-s}$, as before. We formally write $1/S\var(s)=
\sum_n b_nn^{-s}$. The coefficients $b_n$ are given by $b_1=1$,
$$b_n=\sum_{d_1d_2\cdots d_k=n}(-1)^ka_{d_1}a_{d_2}\cdots a_{d_k},\qquad
n>1,$$
where the sum runs over all the finitely many possible product decompositions
of $n$, where the various factors $d_j$ are integers $>1$. The series $\sum_n
b_n n^{-s}$ converges absolutely in some half-plane $\Re s>\sigma_0$, as is
readily seen by considering the extremal situation where $a_n\le0$ for all
$n>1$. In fact, any $\sigma_0$ with $\sum_{n>1}|a_n|n^{-\sigma_0}\le a_1=1$
will do. Multiplying $S\var$ with $1/S\var$, we read off from the coefficients
that $b_1a_1=1$ and
$$\sum_{d:d|n}b_da_{n/d}=0,\qquad n>1.\tag{5-2}$$
Formally, we have $1/S\var=S\psi$, where $\psi(x)=\sum^\infty_{n=1}b_n
e_n(x)$. It follows from \thetag{5-2} that
the functions $\psi_n(x)=\sum_{d:d|n}\bar b_{n/d}e_d(x)$ (they are {\sl not}
the dilates of $\psi$) form a system which is {\sl biorthogonal} to the
original system $\{\var_n\}_n$, with $\var_n(x)=\var(nx)$, that is,
$\langle\var_j,\psi_k\rangle=\delta_{j,k}$, where $\delta_{j,k}$ is
the Kronecker delta symbol, which is $1$ when $j=k$, and $0$ otherwise. It
follows that the system $\{\var_n\}_n$ is {\sl minimal}, that is, each
$\var_n$ lies outside the closure of the linear span of the other vectors.
Moreover, $e_n(x)=\sum_{d:d|n}a_{n/d}\psi_d(x)$, and since these sums contain
only a finite number of terms, it is clear that {\sl the biorthogonal system
$\{\psi_n\}_n$ is complete in} $L^2(0,1)$. The biorthogonal system consists
of dilations of a single function only in the trivial case $\var(x)=e_1(x)=
\sqrt{2}\sin(\pi x)$.

\subhead 5.2. Riesz bases\endsubhead Recall that a basis $\{f_n\}_n$ in a
separable Hilbert space $H$ is a Riesz basis if for some bounded invertible
operator $L$ on $H$, the sequence $\{Lf_n\}_n$ is an orthonormal basis.
An equivalent characterization  is the following (see
\cite{{\bf\Y}, pp. 30--37}).

\proclaim{Lemma 5.1} Let $H$ be a separable Hilbert space. A system
$\{f_n\}_n$ of vectors from $H$  is a
Riesz basis in $H$ if and only if
\roster
\item"{(a)}"  every $f\in H$ can be expanded as $f=\sum_n
c_nf_n$.
\item"{(b)}" there are constants $A$ and $B$, $0<A\leq B<+\infty$, such
that
$$A\Big(\sum_n|c_n|^2\Big)^{1/2}\leq\Big\|\sum_n c_nf_n\Big\|_{L^2}
\leq B\Big(\sum_n |c_n|^2\Big)^{1/2}\tag{5-3}$$
for every finite sequence of scalars $c_{n}$.
\endroster
\endproclaim

We now state the main result of this section.

\proclaim{Theorem 5.2} The system $\var(x)$, $\var(2x)$, $\var(3x)$,$\ldots$
is  a Riesz basis in $L^2(0,1)$ if and only if both $S\var$ and $1/S\var$
belong to $\calM$.
\endproclaim

\demo{Proof} As before, we write $e_{n}(x)=\sqrt{2}\sin( n\pi x)$, and let
$T_{\var}$ be the linear mapping which sends $e_n$ to $\var_n$, for $n=1,\,2,
\,3,\ldots$. Let $\frakF$ denote the dense subspace of $L^2(0,1)$ consisting
of functions $f=\sum_n c_ne_n$, where all but finitely many of the
coefficients $c_n$ equal $0$. By Lemma 5.1, we have a Riesz basis if and only
if the image of $\frakF$ under $T_{\var}$ is dense in $L^2(0,1)$, and
$$A\,\|f\|_{L^2}\le\|T_{\var}f\|_{L^2}\le B\,\|f\|_{L^2},\qquad f\in\frakF,
\tag 5-4$$
holds for some constants $A$ and $B$, $0<A\leq B<+\infty$. After an
application of the transformation $S$, the Riesz basis condition becomes, if
we use \thetag{5-1}, that the image of $S\frakF$ under multiplication by
$S\var$ should be dense in $\calH$, and that
$$A\,\|Sf\|_{\calH}\le\|S\var\,Sf\|_{\calH}\le B\,\|Sf\|_{\calH},\qquad
f\in\frakF.\tag 5-5$$
We first do the sufficiency part. If $S\var$ and $1/S\var$ are both in
$\calM$, then clearly, $S\var\,S\frakF$ is dense, and \thetag{5-5} holds with
$A=\|1/S\var\|_{\calM}^{-1}$ and $B=\|S\var\|_{\calM}$.

We turn to the necessity part. By an argument involving Cauchy sequences,
\thetag{5-4} extends to all $f\in L^2(0,1)$, so that $S\var$ is a multiplier
on $\calH$. By \thetag{5-5} and the fact that $S\var\,S\frakF$ is dense in 
$\calH$, $1/S\var$ is a multiplier on $\calH$ as well. The proof is complete.
\qed
\enddemo

\demo{Example} In conjunction with Theorem 3.1, the theorem above provides a
precise statement about the subtle dependence on the ``closeness'' of $\var$
to $e_1(x)=\sqrt{2}\sin(\pi x)$. To illuminate the condition, consider
$$\var(x)=\sum^\infty_{n=1}\frac{e_n(x)}{n^\tau}=\sqrt{2}\sum^\infty_{n=1}
\frac{\sin(n\pi x)}{n^\tau},$$
for $\tau>\frac12$. Then
$$S\var(s)=\sum^\infty_{n=1}n^{-\tau} n^{-s}=\zeta(\tau+s),$$
and
$$1/S\var(s)=\sum_{n=1}^\infty\mu(n)n^{-\tau} n^{-s},$$
where $\mu(n)$ is the M\"{o}bius function, which has $\mu(1)=1$, $\mu(n)=
(-1)^k$ when $n$ is the product of $k$ different primes, and $\mu(n)=0$ if $n$
is divisible by a square (other than $1$). The function $S\var(s)$ is bounded
in the half-plane $\Re s>0$ if and only if $\tau>1$. Using Theorem 5.2 and the
fact that multipliers are bounded analytic functions, we conclude that the
system $\var(x),\var(2x),\ldots$ is a Riesz basis in $L^2(0,1)$ if and only
if $\tau>1$.
\enddemo

The above example is covered by the following corollary to Theorem 5.2.

\proclaim{Corollary 5.3} If the coefficients $a_n$ of $S\var\in\calH$ are
totally multiplicative, the functions $\var(x),\,\var(2x),\,\var(3x),\ldots$
form a Riesz basis in $L^2(0,1)$ if and only if $\sum_p|a_p|<+\infty$, where
the sum runs over the primes.
\endproclaim

\demo{Proof} Since the totally multiplicative coefficients come from a
function in $\calH$, they must satisfy $\sup_n|a_n|<1$, and by the Euler
product formula,
$$S\var(s)=\sum_{n=1}^\infty a_nn^{-s}=\prod_{p}\big(1-a_pp^{-s}\big)^{-1}.
\tag 5-6$$
By the remark following the proof of the multiplier theorem, $\calM$ is
isometrically isomorphic to $H^\infty(\D^\infty)$. As we calculate the norm
of $\fQ S\var$ there, using the analog of \thetag{5-6}, we obtain
$$\|S\var\|_{\calM}=\prod_{p}\big(1-|a_p|\big)^{-1},$$
so that in view of Theorem 5.2, $\sum_p|a_p|<+\infty$ is certainly a necessary
condition to have a Riesz basis. Incidentally, this is also
sufficient, since
$$
 1/ S\var(s)=\sum_{n=1}^{\infty}\mu(n) a_{n} n^{-s}=
  \prod_{p}\big(1-a_{p}p^{-s}\big).\ \text{\qed}
$$
% \|_{\infty} \leq  On the other hand, if we take
% logarithms on both sides of \thetag{5-6}, and compute norms in $\calM$, then
% it follows that $\sum_p|a_p|<+\infty$ is sufficient for $\log S\var$ to be in
% $\calM$. Since $\calM$ is a Banach algebra, $S\var$ and $1/S\var$ are then
% both in $\calM$. The proof is complete.
 \enddemo

Note next the following consequence of Theorems 3.1 and 5.2; for the second
statement, one should keep in mind the identification of $\calM$ with $H^\infty
(\D^\infty)$ supplied by Section 3.

\proclaim{Corollary 5.4} If $\var(x)=\sum_n a_ne_n(x)$, and $\sum_{n>1}|a_n|
<a_1=1$, then the collection $\{\var(nx)\}_n$ $(n=1,2,3,\ldots)$ forms a Riesz
basis in $L^2(0,1)$. On the other hand, if $a_n=0$ unless $n$ is a prime or
$1$, the condition $\sum_{n=2}^{\infty}|a_n|<1$ is necessary to have a Riesz
basis.
\endproclaim

\subhead 5.3. The  completeness problem\endsubhead The completeness problem
is considerably more delicate than the Riesz basis problem. It does not seem
likely that it has a solution as simple as the one given by Theorem 5.2.
In fact, as we shall see, it is equivalent to the problem of describing the
cyclic vectors in the space $H^2(\D^\infty)$, which is definitely quite hard.
The relationship to the invariant subspaces in an infinite-dimensional setting
was pointed out by the editors of Beurling's collected works \cite{\bf\Be}.
The level of difficulty of Beurling's completeness problem is emphasized by the
fact that a complete characterization of the cyclic vectors in $H^2(\D^N)$ is
known only when the (complex) dimension of the polydisk is $N=1$.

A subspace $I$ of $\calH$ is said to be {\sl invariant} if it is closed, and
$fg\in I$ whenever $f\in I$ and $g\in\frakG$. One shows, by suitably
approximating functions in $\calM$ with elements of $\frakG$, that invariant
subspaces are actually invariant under multiplication by elements of $\calM$.
The invariant subspace generated by a function $f\in\calH$ is denoted by $[f]
_{\calH}$, and we say that $f$ is cyclic if $[f]_{\calH}=\calH$.
Similarly, a subspace $J$ of $H^2(\D^\infty)$ is called {\sl invariant} if it
is closed, and invariant under multiplication by polynomials in the coordinate
functions $z_1,\,z_2,\ldots$ on $\D^\infty$. The invariant subspaces are
actually invariant under multiplication by $H^\infty(\D^\infty)$. We write
$[f]_{H^2(\D^\infty)}$ for the invariant subspace generated by $f\in
H^2(\D^\infty)$, and say that $f$ is cyclic if $[f]_{H^2(\D^\infty)}=
H^2(\D^\infty)$.

A moment's reflection on the definition of $T_\var$ reveals that $\var_1,\,
\var_2,\,\var_3,\ldots$ form a complete system in $L^2(0,1)$ if and only if
the image of $\frakF$ under $T_\var$ is dense in $L^2(0,1)$. After an
application of the transformation $S$, using \thetag{5-1}, the latter condition
becomes the requirement that $S\var$ times the set $\frakG=S\frakF$ of finite
Dirichlet series be dense in $\calH$, {\sl id est}, that $S\var$ be cyclic in
$\calH$. By the isometry between $\calH$ and $H^2(\D^\infty)$ supplied by the
operator $\fQ$ back in Section 2, the following can be said.

\proclaim{Theorem 5.5} Let $\var$ be as above. Then the following are
equivalent:
\roster
\item"{(a)}" the system $\var_1,\,\var_2,\,\var_3,\ldots$ is complete in
$L^2(0,1)$.
\item"{(b)}" the function $S\var$ is cyclic in $\calH$, that is, the subspace
$S\var\,\calM$ is dense in $\calH$.
\item"{(c)}" the function $\fQ S\var$ is cyclic in $H^2(\D^\infty)$.
\endroster
\endproclaim

A necessary condition for the completeness of the system $\{\var_n\}_n$
mentioned by Beurling in his 1945 seminar ({\sl confer} \cite{{\bf\Be}, p.
378}), is that the function $S\var(s)$ be zero-free in the half-plane $\Re s>
\frac12$. We shall prove something slightly stronger. For a sequence
$z=(z_j)_j$ in $\D^\infty$, let $\|z\|_{l^\infty}=\sup\{|z_j|:\,j=1,2,\ldots
\}$ and $\|z\|^2_{l^2}=\sum_j|z_j|^2$. We extend this notation to sets
$\Omega\in\D^\infty$ by taking suprema over all the points in $\Omega$.

\proclaim{Lemma 5.6} If the system $\var_1,\var_2,\ldots$ is complete in
$L^2(0,1)$, and $\Omega\subset\D^{\infty}$ has $\|\Omega\|_{l^\infty}<1$ and
$\|\Omega\|_{l^2}<+\infty$, then
$$\inf_{z\in\Omega}|\fQ S\var(z)|>0.$$
\endproclaim

\demo{Proof} In view of Theorem 5.5, we need to show that if $g\in
H^2(\D^\infty)$ is cyclic, then $|g|$ is bounded away from $0$ on $\Omega$.
We argue by contradiction. So, suppose $|g|$ is not bounded away from $0$ on
$\Omega$; then there is a sequence of points $z(k)$ in $\Omega$ such that
$g(z(k))\to0$ as $k\to+\infty$. When $l^2(\N)$ is given the weak topology,
its closed unit ball becomes a compact metric space (the Banach-Alaoglu
theorem), so that by the Bolzano-Weierstra\ss{} theorem, the sequence
$\{z(k)\}_k$ possesses a cluster point $z(\infty)$, with
$\|z(\infty)\|_{l^\infty}<1$ and $\|z(\infty)\|_{l^2}<+\infty$. By the
conditions on $\Omega$, $\sup_{z\in\Omega}\prod_j(1-|z_j|^2)^{-1}<+\infty$,
so that point evaluations in $\Omega$ are uniformly bounded. It follows that
as the subsequence of points $z(k_l)$ converges to $z(\infty)$, the function
values $g(z(k_l))$ converge to $g(z(\infty))$, so that $g(z(\infty))=0$.
We conclude that $g$ cannot be cyclic, as it is annihilated by a bounded
point evaluation.
\qed
\enddemo

\demo{Remark} In particular, it follows that for the system $\var_1,\var_2,
\ldots$ to be complete, a necessary condition is that
$$\sup_{\Re s>\sigma_0}|S\var(s)|>0,$$
for any $\sigma_0>\frac12$. This result is sharp, in the sense that for
complete systems, $\inf_{\Re s>\frac12}|S\var(s)|$ $=0$ can occur. An example
is given after Corollary 5.8. A slightly more refined necessary condition for
completeness can be obtained from Theorem 4.10: for all characters $\chi$,
the function $(S\var)_\chi(s)=\sum_n a_n\chi(n)n^{-s}$ must be outer in
$H^2_\infty(\C_{\frac12})$. Also, by a remark Helson makes in \cite{\bf\Heb},
it is necessary that $(S\var)_\chi$ almost surely be outer in $H^2_{\text{i}}
(\C_+)$.
\enddemo

\subhead 5.4. The Dirichlet-type space $\calH_d$\endsubhead
We introduce next a subspace of $\calH$ which seems to play a natural role
in the study of the completeness problem. To this end, let
$$d(n)=\sum_{d:d|n}1$$
denote the number of divisors of a natural number $n$. As a function of $n$,
the number of divisors is quite irregular, but it is well-known that $d(n)=
O(n^\delta)$ holds for every $\delta>0$; however, the average order of $d(n)$
is $\log n$, {\sl confer} \cite{{\bf\HW}, pp. 260--266}. We denote by
$\calH_d$ the collection of Dirichlet series $f(s)=\sum^\infty_{n=1}a_nn^{-s}$
for which
$$\sum^\infty_{n=1} |a_n|^2d(n)<\infty.$$
The set of those $f$ in $\calH_d$ for which $a_n=0$ unless $n$ is a power of 
a given prime $p$, coincides with the usual Dirichlet space on the unit disk,
if we set $z=p^{-s}$.

Clearly, $\calH_d\subset\calH$. For the interesting class of series with
totally multiplicative coefficients ($a_{mn}=a_ma_n$), we have the opposite
inclusion. Indeed,
$$\sum_n|a_n|^2d(n)=\left(\sum_n|a_n|^2\right)^2$$
in this favorable case. The identity can be verified by expanding both sums
as Euler products. If we write $n=\prod_p p^{\nu_p}$, where the product runs
over the primes, and all but finitely many of the nonnegative integers $\nu_p$
are $0$, we have
$$\multline
\sum_n|a_n|^2d(n)=\sum_n|a^{\nu_2}_2a^{\nu_3}_3\cdots a^{\nu_p}_p
\cdots|^2(1+\nu_2)(1+\nu_3)\cdots(1+\nu_p)\cdots\\
=\sum^\infty_{\nu_2=0}(1+\nu_2)|a_2|^{2\nu_2}\cdot
\sum^\infty_{\nu_3=0}(1+\nu_3)|a_3|^{2\nu_3}\cdots
\sum^\infty_{\nu_p=0}(1+\nu_p)|a_p|^{2\nu_p}\cdots\\
=\big(1-|a_2|^2\big)^{-2}\cdot\big(1-|a_3|^2\big)^{-2}\cdots
\big(1-|a_p|^2\big)^{-2}\cdots=\prod_p\big(1-|a_p|^2\big)^{-2},
\endmultline$$
provided that each $|a_p|<1$, and $a_1=1$. Note that if some $|a_p|\geq 1$,
the sequence $\{a_n\}_n$ is not even in $l^2(\N)$. We have, by using a similar
argument,
$$\sum^\infty_{n=1}|a_n|^2=\prod_p\big(1-|a_p|^2\big)^{-1},$$
where $a_1=1$ and $|a_p|<1$. The product converges if and only if
$$\sum_p|a_p|^2<+\infty, \ \ |a_p|<1.$$
{\sl This is a necessary and sufficient criterion for a Dirichlet series with
totally multiplicative coefficients to belong to the space $\calH$, and thus
to the space $\calH_{d}$}.

The reciprocal series of a Dirichlet series $f(s)=\sum_{n}a_nn^{-s}$ with
totally multiplicative coefficients and $a_1=1$ is
$$\frac1{f(s)}=\sum^\infty_{n=1}\mu(n)a_nn^{-s},$$
where $\mu(n)$ is the M\"{o}bius function. A direct calculation yields
$$\sum^\infty_{n=1}|\mu(n)a_n|^2=\prod_p\left(1+|a_p|^2\right),\qquad
\sum_n|\mu(n)a_n|^2d(n)=\prod_p\left(1+2|a_p|^2\right) $$
for the expressions defining the norms in $\calH$ and in $\calH_d$.
These norms are finite if and only if $\sum_p|a_p|^2<+\infty$. This criterion
for $1/f(s)$ differs from the one for $f(s)$ in one respect: the condition
$|a_p|<1$ is skipped.

\subhead 5.5 Completeness, continued\endsubhead We now use the space
$\calH_{d}$ to formulate a fairly general sufficient condition for
completeness.

\proclaim{Theorem 5.7} The system $\{\var_n\}_n$ is complete in $L^2(0,1)$ if
one of the following conditions is fulfilled:
\roster
\item"{(a)}" $S\var$ and $1/S\var$ are both in $\calH_d$.
\item"{(b)}" $S\var\in\calM$ and $1/S\var\in\calH$.
\item"{(c)}" $S\var\in\calH$ and $1/S\var\in\calM$.
\endroster
\endproclaim

\demo{Proof} We first look at parts (a) and (b). As before, we write
$1/S\var(s)=\sum_n b_n n^{-s}$. Put $R_N(s)=\sum_{n=1}^Nb_nn^{-s}$.
Suppose we can show that $\|S\var R_N\|_{\calH}$ has a bound which is
independent of $N$. Then $S\var R_N$ tends to the constant function $1$
weakly, as $N\to+\infty$. By standard functional analysis, $1$ is in the
norm closure of $S\var\,\calM$, so that $S\var$ is cyclic in $\calH$. Under
condition (b), it is clear that this boundedness condition is fulfilled,
by the estimate $\|R_N\|_{\calH}\le\|1/S\var\|_{\calH}$. Under condition
(a), we apply the Cauchy-Schwarz inequality,
$$\multline
\|S\var R_N\|^2_{\calH}=\sum_k\Big|\sum_{m,n:\atop mn=k,m\le N}a_nb_m\Big|^2
\le\sum_{k}d(k)\sum_{m,n:mn=k}|a_n|^2|b_m|^2\\
\leq\left(\sum_{n}|a_n|^2d(n)\right)\left(\sum_{m}|b_m|^2d(m)\right)
=\|S\var\|_{\calH}^2\|1/S\var\|_{\calH}^2,
\endmultline$$
where we used that $d(mn)\le d(m)d(n)$. Part (a) follows. Under condition (c),
we just multiply $S\var$ with the function $1/S\var$ in $\calM$, to get that
$1\in[S\var]_{\calH}$, so that $S\var$ is cyclic. The proof is complete.
\qed
\enddemo

\proclaim{Corollary 5.8} If the coefficients of $S\var\in\calH$ are totally
multiplicative, the system $\{\var_n\}_n$ is complete in $L^2(0,1)$.
\endproclaim

\demo{Proof} According to the previous subsection, we have that both $S\var$
and $1/S\var$ are in $\calH_d$, because the coefficients are totally
multiplicative. The assertion now follows from part (a) of Theorem 5.7.
\qed
\enddemo

\demo{Remark} It follows that in the context of the example considered after
Theorem 5.2, we have completeness if and only if $\tau>\frac12$. In a way,
this reproves a theorem of Helson \cite{\bf\Heb}: the function $S\var$ is
cyclic in $\calH$ if $\tau>\frac12$.
\enddemo

\demo{Example} We use Theorem 5.7 to construct a complete system $\var_1,
\var_2,\ldots$ with
$$\inf_{\Re s>\frac12}\big|S\var(s)\big|=0,$$
as promised in the remark following Lemma 5.6. Let $b_1=1$, $b_2=0$, and put
$b_p=p^{-1/2}(\log\log p)^{2/3}$ for prime indices $p>2$, and extend the
sequence multiplicatively, so that $b_{mn} =b_mb_n$. By Chebyshev's theorem
\cite{{\bf\La}, Vol. 1, p. 25}, the series
$$\sum_{p}|b_p|^2=\sum_{p:p>2}\frac1{p\,(\log\log p)^{4/3}}$$
(summation over the primes) converges if and only if
$$\sum^\infty_{n=3}\frac1{n(\log\log n)^{4/3}\log n}$$
converges. Since the last series converges by the integral test, we conclude
from the results of the previous subsection that the functions $\Phi$ and
$1/\Phi$ are in $\calH_d$, where $\Phi(s)=\sum_nb_nn^{-s}$. Let $\var$ be given
by $S\var=1/\Phi$. The completeness now follows from Theorem 5.7. Again by
Chebyshev's theorem,
$$\Phi(\sigma)\geq\sum_{p:p>2}\frac1{p^{1/2+\sigma}(\log\log p)^{2/3}}\ge A
\sum_{n=3}^\infty\frac1{n^{1/2+\sigma}(\log\log n)^{2/3}\log n},$$
for some positive absolute constant $A$, whence $\Phi(\sigma)\to+\infty$ as 
$\sigma\to\frac12$, since the right hand side series diverges for $\sigma=
\frac12$. It follows that $S\var(\sigma)\to0$ as $\sigma\to\frac12$.
\enddemo

\proclaim{Theorem 5.9} Let $S\var(s)=1+\sum_pa_pp^{-s}$, summing only over
the primes. Then the system $\{\var_n\}_n$ is complete in $L^2(0,1)$ if and
only if $\sum_p|a_p|\le1$.
\endproclaim

\demo{Proof} To see the necessity, observe that if $\sum_{p}|a_{p}|>1$, there
exist a positive integer $N$ and and a point $z=(z_1,z_2,\ldots)\in\D^\infty$
with $z_j=0$ for $j>N$, such that ($p_j$ is the $j$-th prime)
$$\fQ S\var(z)=1+\sum_{j=1}^{\infty}a_{p_{j}}z_{j}=0.$$
By Lemma 5.6, we cannot have completeness.

To get the sufficiency, we introduce the auxiliary function $g=1-S\var$, and
note that $\|g\|_{\calM}\le1$. The function $1-g^n=(1-g)(1+g+\ldots+g^{n-1})$
tends to $1$ in the norm of $\calH$ as $n\to+\infty$, and for each $n$, it
is in $S\var\,\calM$. It follows that $S\var$ is cyclic.
\qed \enddemo

\head
6. Possible directions of further investigation
\endhead

According to a theorem of Helson (Theorem 4.4), we have that for $f\in\calH$, 
with series expansion \thetag{2-1}, the Dirichlet series
$$f_\chi(s)=\sum_{n=1}^\infty a_n\chi(n)n^{-s}$$
converges in the half-plane $\Re s>0$, for almost all characters $\chi$. In
terms of the coefficients, this  amounts to having, almost surely in $\chi$,
$$\sum_{n=1}^N a_n\chi(n)=O\big(N^{\vare}\big)\qquad\text{as}\quad
N\to+\infty,\tag 6-1$$
for each fixed $\vare>0$. In view of Theorem 4.5 and Lennart Carleson's 
theorem on the almost everywhere convergence of Fourier series with $l^2$ 
coefficients \cite{\bf\Carle}, which cover the probabilistic and deterministic
extremes, the following conjecture seems reasonable: for almost all $\chi$,
$$\sum_{n=1}^N a_n\chi(n)=O\big(1\big)\qquad\text{as}\quad N\to+\infty.$$
The conjecture appears to be related to the question of how many of the 
stochastic variables $\chi(n)$ on a given interval $N\le n<N+k$ are mutually 
independent. The latter translates into the following basic number-theoretic 
problem: What is the multiplicatively rational dimension of the interval 
$N\le n<N+k$? Here, the multiplicatively rational dimension of a set $E$ of 
positive integers is the dimension of the $\Q$-linear span of $\log E$ over 
the field $\Q$. For instance, the multiplicatively rational dimension of the 
interval $1\le n\le k$ equals the number of primes $\le k$. 

By Corollary 4.7, almost every $\zeta_\chi(s)$ converges to a zero-free 
analytic function in the half-plane $\Re s>\frac12$. It will now be indicated 
how this result may be used to obtain information about particular characters.
Given positive reals $M$ and $\vare$, let $\Omega(\vare,M)$ be the set of all 
$\chi$ for which 
$$\Big|\sum_{n=1}^{N}\chi(n)\Big|\le M\,\cdot N^{1/2+\vare},$$
for all $N=1,2,3,\ldots$. By Corollary 4.8, the $\rho$-mass of $\Omega(\vare,
M)$ tends to $1$ as $M\to+\infty$, for fixed $\vare$. A point $\chi_0$ in 
$\Char$ is a {\sl mass point for} $\Omega(\vare,M)$ provided that each open 
neighborhood of $\chi_0$, intersected with $\Omega(\vare,M)$, has positive 
$\rho$-mass. {\sl If $\chi_0$ is a mass point for some $\Omega(\vare,M)$, then
$\zeta_{\chi_0}$ extends analytically to a zero-free function on} $\Re s>
\frac12+\vare$. This is so because this function may be approximated by 
functions $\zeta_\chi(s)$ which lack zeros on $\Re s>\frac12$, in the topology
of uniform convergence on compact subsets of $\Re s>\frac12+\vare$. This may 
be a way to handle non-principal Dirichlet characters, a topic to be 
discussed below.

As a matter of definition, a character is a multiplicative mapping $\chi:
\Q_+\to\T$, and if we like, we may extend it to all of $\Q\setminus\{0\}$ by
setting $\chi(-1)$ equal to $1$ or $-1$. The characters that so extend
continuously to the archimedian completion $\R\setminus\{0\}$ are of the type
$\chi(r)=|r|^{-it}$ for some real number $t$ (if we fix $\chi(-1)=1$). There
are also the non-archimedian completions $\Q_p$, the $p$-adic number field,
for a given fixed prime $p$. The {\sl Dirichlet characters} associated with
$\Q_p$ are the ones that extend continuously to $\Q_p\setminus\{0\}$, and one
checks that for some $k\in\N$, when restricted to the integers, they are
periodic with period $p^k$, except on numbers divisible by $p$. Traditionally,
one sets the value of the Dirichlet character $\chi$ equal to $0$ on the
integers divisible by $p$, but this is not our choice here. The principal
Dirichlet characters take value $1$ on all integers indivisible by $p$ (there
are several because we do not prescribe the value at $p$). The nonprincipal
Dirichlet characters all have partial sums $\sum_{n:n\le N}\chi(n)$ that grow
like $O(\log N)$ as $N$ tends to infinity (with the traditional definition,
these sums are actually bounded), so at first glance they seem like ideal
candidates for mass points. In the final analysis, this may prove to be
na\"\i ve.

\subhead Acknowledgement\endsubhead The first author wishes to thank Paul
Cohen (Stanford) for helpful conversations at an early stage of this
work.

\enddocument